\documentclass[10pt]{article}

\usepackage{amsmath,amsxtra,amssymb,latexsym, amscd,amsthm}
\usepackage{here}
\usepackage{graphpap}
\usepackage{graphicx}
\usepackage{color}
\usepackage[]{algorithm2e}

\newtheorem{theorem}{Theorem}[section]

\newtheorem{lemma}[theorem]{Lemma}

\theoremstyle{definition}

\newtheorem{remark}[theorem]{Remark}

\numberwithin{equation}{section}

\def\<{\langle}
\def\>{\rangle}

\def\limn{\lim_{n\to\infty}}
\def\limsupn{\limsup_{n\to\infty}}
\def\liminfn{\liminf_{n\to\infty}}
\def\={&=&}
\definecolor{purple}{rgb}{0.4, 0.0, 0.4}

\begin{document}

\begin{center}
\large \bf Variational method for multiple parameter identification in elliptic PDEs
\end{center}

\vspace{0.1cm}

\centerline{Tran Nhan Tam Quyen\let\thefootnote\relax\footnote{Email: quyen.tran@uni-hamburg.de,~ quyen.tran@uni-goettingen.de}\let\thefootnote\relax\footnote{The author gratefully acknowledges support of the Alexander von Humboldt-Foundation and the Lothar Collatz Center for Computing in Science at the University of Hamburg, and the University of Goettingen, State of Lower Saxony, Germany}
}

Department of Mathematics, University of Hamburg, Bundesstra\ss{}e 55, D-20146 Hamburg, Germany\\
Institut f\"ur Numerische und Angewandte Mathematik, Universit\"at G\"ottingen, Lotzestra\ss{}e 16-18, D-37083 G\"ottingen, Germany

\centerline{\rule{8cm}{.1pt}}

\vspace{0.3cm}

{\small {\bf Abstract:} In the present paper we investigate the inverse problem of
identifying simultaneously the diffusion matrix, source term and boundary condition in the Neumann boundary value problem for an elliptic partial differential equation (PDE) from a measurement data, which is weaker than required of the exact state. A variational method based on  energy functions with Tikhonov regularization is here proposed to treat the identification problem. We discretize the PDE with the finite element method and prove the convergence as well as analyse error bounds of this approach.}



\section{Introduction}

Let $\Omega$ be an open bounded connected domain of $R^d,
1\le d \le 3$ with polygonal boundary $\partial \Omega$. In this paper we study the problem of identifying simultaneously the {\it diffusion matrix $Q$, source term $f$} and {\it boundary condition $g$} as well as the {\it state} $\Phi$
in the Neumann boundary value problem for the elliptic PDE
\begin{eqnarray}
-\nabla \cdot (Q\nabla \Phi) &=& f \quad\mbox{in}\quad \Omega,  \label{ct1}\\
Q\nabla \Phi\cdot\vec{n} &=& g \quad\mbox{on}\quad \partial\Omega \label{ct2}
\end{eqnarray}
from a measurement $z_\delta\in L^2(\Omega)$ of the solution $\Phi\in H^1(\Omega)$, where $\vec{n}$ is the unit outward normal on $\partial\Omega$.

To formulate precisely our problem, let us first denote by $\mathcal{S}_d$ the set of all symmetric, real $d \times
d$-matrices equipped with the inner product $M \cdot N :=
\mbox{trace} (MN)$ and the corresponding norm $\| M \|_{\mathcal{S}_d} = (M \cdot
M)^{1/2} = \left( \sum_{i,j = 1}^d m_{ij}^2\right)^{1/2}$, where $M:= (m_{ij})_{i,j=\overline{1,d}}$.
Furthermore, for $1 \le p \le \infty$ we set
\begin{eqnarray*}
\mathbf{L}^p_{ \mbox{\tiny sym}}(\Omega) &:=& \left\{H := (h_{ij})_{i,j=\overline{1,d}}\in {L^p(\Omega)}^{d \times d} \quad\big|\quad H(x):= (h_{ij}(x))_{i,j=\overline{1,d}} \in \mathcal{S}_d \quad\mbox{a.e. in}\quad \Omega \right\}.
\end{eqnarray*}
In $\mathbf{L}^2_{ \mbox{\tiny sym}}(\Omega)$ we use the scalar product
$\left ( H^1, H^2 \right)_{\mathbf{L}^2_{ \mbox{\tiny sym}}(\Omega)} = \sum_{i,j = 1}^d  ( h^1_{ ij}, h^2_{ ij} )_{L^2(\Omega)}$
and the corresponding norm
$\|H\|_{\mathbf{L}^2_{ \mbox{\tiny sym}}(\Omega)} := \left( \sum_{i,j = 1}^d \| h_{ ij}\|^2_{L^2(\Omega)}\right)^{1/2} = \left( \int_\Omega \| H(x) \|^2_{\mathcal{S}_d}dx\right)^{1/2}$, while the space $\mathbf{L}^\infty_{ \mbox{\tiny sym}}(\Omega)$ is endowed with the norm $\|H\|_{\mathbf{L}^\infty_{ \mbox{\tiny sym}}(\Omega)} := \max_{i,j=\overline{1,d}}
\|h_{ij}\|_{L^{\infty}(\Omega)}$.

Let us denote by
\begin{eqnarray*}
\mathcal{H}_{ad} &:=& \mathcal{Q}_{ad} \times \mathcal{F}_{ad} \times \mathcal{G}_{ad}
\end{eqnarray*}
with
\begin{eqnarray}\label{5/12/12:ct3}
\mathcal{Q}_{ad} &:=& \left\{ Q \in \mathbf{L}^\infty_{ \mbox{\tiny sym}}(\Omega)
\quad\big|\quad \underline{q} |\xi|^2
\le Q(x) \xi \cdot \xi \le \overline{q} |\xi|^2 \quad\mbox{for all}\quad \xi \in R^d \right\}, \\
\mathcal{F}_{ad} &:=& L^2(\Omega),\nonumber\\
\mathcal{G}_{ad} &:=& L^2(\partial\Omega)\nonumber
\end{eqnarray}
and $\underline{q}$, $\overline{q}$ being given constants satisfying $\overline{q} \ge \underline{q} >0$. Let 
$$\gamma :\quad H^1(\Omega) \quad\to\quad H^{1/2}(\partial\Omega)$$
be the continuous Dirichlet trace operator and $ H^1_\diamond(\Omega)$ be the closed subspace of $H^1(\Omega)$ consisting all functions with zero-mean on the boundary, i.e.
\begin{eqnarray*}
 H^1_\diamond(\Omega) &:=& \left\{ u \in H^1(\Omega) \quad\Big|\quad \int_{\partial\Omega} \gamma udx =0\right\}
\end{eqnarray*}
while $C_\Omega$ stands for the positive constant appearing in the Poincar\'e-Friedrichs inequality (cf.\  \cite{Pechstein})
\begin{eqnarray}\label{20-10-15ct1}
C_\Omega \int_\Omega \varphi^2dx &\le& \int_\Omega |\nabla \varphi|^2dx \quad\mbox{for all}\quad \varphi\in  H^1_\diamond(\Omega).
\end{eqnarray}
Then, due to the coervicity condition
\begin{eqnarray}\label{coercivity}
\| \varphi\|^2_{H^1(\Omega)} &\le& \frac{1+ C_\Omega}{C_\Omega} \int_\Omega | \nabla \varphi|^2dx \quad\le\quad \frac{1+ C_\Omega}{C_\Omega \underline{q}} \int_\Omega Q\nabla\varphi \cdot \nabla\varphi dx
\end{eqnarray}
holding for all $\varphi\in  H^1_\diamond(\Omega), Q\in\mathcal{Q}_{ad}$
and the Lax-Milgram lemma, we conclude for each $(Q,f,g) \in \mathcal{H}_{ad}$,
there exists a unique weak solution $\Phi$ of (\ref{ct1})--(\ref{ct2}) in the sense that $\Phi\in H^1_\diamond(\Omega)$ and satisfies the identity
\begin{eqnarray}
\int_\Omega  Q \nabla \Phi \cdot \nabla \varphi dx &=&
(f,\varphi) + \langle g,\gamma\varphi\rangle \label{ct9}
\end{eqnarray}
for all $\varphi\in H^1_\diamond(\Omega)$. Here
the expressions $(\cdot,\cdot)$ and $\langle\cdot,\cdot\rangle$ stand for the scalar product on space $L^2(\Omega)$ and $L^2(\partial\Omega)$, respectively. Furthermore, there holds the priori estimate
\begin{eqnarray}\label{mq5}
\left\|\Phi\right\|_{H^1(\Omega)}
&\le& \frac{1+ C_\Omega}{C_\Omega \underline{q}} \left( \left\|\gamma\right\|_{\mathcal{L}\big(H^1(\Omega),H^{1/2}(\partial\Omega)\big)} \left\|g\right\|_{L^2(\partial\Omega)} + \|f\|_{L^2(\Omega)}\right) \nonumber\\
&\le& C_{\mathcal{N}} \left( \left\|g\right\|_{L^2(\partial\Omega)} + \|f\|_{L^2(\Omega)}\right)
\end{eqnarray}
with
\begin{eqnarray*} C_{\mathcal{N}} &:=& \frac{1+ C_\Omega}{C_\Omega \underline{q}} \max \left( 1, \left\|\gamma\right\|_{\mathcal{L}\big(H^1(\Omega),H^{1/2}(\partial\Omega)\big)}\right).
\end{eqnarray*}
Then we can define the {\it non-linear coefficient-to-solution operator}
\begin{eqnarray*}
\mathcal{U} : \quad\mathcal{H}_{ad}
&\rightarrow& H^1_\diamond(\Omega)
\end{eqnarray*}
which maps each $(Q,f,g) \in \mathcal{H}_{ad}$ to the unique weak solution $\mathcal{U}_{Q,f,g} := \Phi$ of the problem (\ref{ct1})--(\ref{ct2}). Here, for convenience in computing numerical solutions of the pure Neumann problem we normalize the solution with vanishing mean on the boundary (cf.,\ e.g., \cite[Subsection 5.2]{hkq}, \cite[Section 2]{bangti2}); however, all results performed in the present paper are still
valid for the normalization of solutions of the Neumann problem with zero-mean over the domain, i.e.
$
\mathcal{U}_{Q,f,g} \in  \big\{ u \in H^1(\Omega) ~\Big|~ \int_{\Omega} udx =0\big\}.
$
The identification problem is now stated as follows:
\begin{center}
{\it Given $\Phi^\dag := \mathcal{U}_{Q,f,g} \in H^1_\diamond(\Omega)$, find an element $(Q,f,g) \in \mathcal{H}_{ad}$\\
such that (\ref{ct9}) is satisfied with $\Phi^\dag$  and $Q, f, g$.}
\end{center}
This inverse problem may have more than one solution and it is highly ill-posed. In fact, assume that the exact $\Phi^\dag \in C^2_c(\Omega)$, the space of all functions having second-order derivatives with compact support in $\Omega$. Then, for all $Q\in {C^1(\Omega)}^{d\times d} \cap \mathcal{Q}_{ad}$ the element $(\bar{Q},\bar{f},\bar{g}) := \left( Q, -\nabla (Q \cdot \nabla \Phi^\dag), 0 \right) $ is a solution of the above identification problem, i.e. $\mathcal{U}_{\bar{Q},\bar{f},\bar{g}} = \Phi^\dag$. In other words we are considering to solve an equation $\mathcal{U}_{Q,f,g} = \Phi^\dag$, where the forward operator $\mathcal{U}$ is non-linear and {\it non-injective}.  Without using
additional objective a-priori information or
without exploiting other observation data as considering here, it is difficult for us to classify sought targets. Following the general convergence theory for ill-posed problems (see, e.g., \cite[Chapter 5]{chavent2009} and \cite[Subsection 3.2.1]{schuster}, or the classical monograph \cite[Section 10.1]{Engl_Hanke_Neubauer}), in the present paper we are interested in finding {\it exact solutions with penalty minimizing}, which is defined as
\begin{eqnarray}\label{8-1-16ct2}
\left( Q^\dag,f^\dag,g^\dag\right) &:=& \arg \min_{(Q, f, g) \in \mathcal{I} \left( \Phi^\dag \right) } \mathcal{R}(Q, f, g),
\end{eqnarray}
where
$
\mathcal{I}(\Phi^\dag) := \left\{ (Q,f,g) \in
\mathcal{H}_{ad} ~|~\mathcal{U}_{Q,f,g} = \Phi^\dag\right\}
$
and the penalty term
\begin{eqnarray*}
\mathcal{R}(Q,f,g) &:=& \|Q\|^2_{\mathbf{L}^2_{ \mbox{\tiny sym}}(\Omega)} + \|f\|^2_{L^2(\Omega)} + \|g\|^2_{L^2(\partial\Omega)}.
\end{eqnarray*}
We note that the admissible set
$\mathcal{I}(\Phi^\dag)$ of the problem (\ref{8-1-16ct2})
is non-empty, convex and weakly closed in $\mathbf{L}^2_{ \mbox{\tiny sym}}(\Omega) \times L^2(\Omega) \times L^2(\partial\Omega)$, so that the minimizer $(Q^\dag, f^\dag, g^\dag)$ is defined uniquely. Furthermore, the exact data $\Phi^\dag$ may not be known in practice, thus we assume instead of $\Phi^\dag$ to have a measurement $z_\delta \in L^2(\Omega)$
such that
\begin{eqnarray}\label{16-9-16ct1}
\left\|\Phi^\dag - z_\delta \right\|_{L^2(\Omega)} &\le& \delta \quad\mbox{with}\quad \delta >0.
\end{eqnarray}
Our identification problem is now to reconstruct $\left( Q^\dag,f^\dag,g^\dag\right) \in \mathcal{H}_{ad}$ from $z_\delta$. 

Let $\left(\mathcal{T}^h\right)_{0<h<1}$ denote a family of triangulations of the domain $\overline{\Omega}$ with the mesh size $h$ and $\mathcal{U}^h$ be the approximation of the operator $\mathcal{U}$ on the piecewise linear, continuous finite element space associated with $\mathcal{T}^h$. Furthermore, let $\Pi^h$ be the Cl\'ement's mollification interpolation operator (cf.\ \S \ref{Finite element method}). The standard method for solving the above mentioned identification problem  is the output least squares one with Tikhonov regularization, i.e. one considers a minimizer of the problem
\begin{eqnarray}\label{30-9-16ct1}
\min_{(Q,f,g) \in \mathcal{H}_{ad}}
\left\|\mathcal{U}^h_{Q,f,g}-
\Pi^h z_\delta\right\|^2_{L^2(\Omega)} + \rho\mathcal{R}(Q,f,g)
\end{eqnarray}
as a discrete approximation of the identified coefficient $\left( Q^\dag,f^\dag,g^\dag\right)$, here $\rho >0$ is the regularization parameter. However, due to the non-linearity of the coefficient-to-solution operator, we are faced with certain difficulties in holding the {\it non-convex} minimization problem (\ref{30-9-16ct1}). Thus, instead of working with the above least squares functional and following the use of energy functions 
(cf.\ \cite{Kohn_Vogelius1,know-wall,wang_zou}), in the present work the {\it convex} cost function (cf.\ \S\ref{Finite element method})
\begin{eqnarray*}
(Q,f,g) \in \mathcal{H}_{ad} &\mapsto& \mathcal{J}^h_\delta (Q,f,g) := \int_\Omega Q\nabla \left(\mathcal{U}^h_{Q,f,g}-
\Pi^h z_\delta\right) \cdot  \nabla \left(\mathcal{U}^h_{Q,f,g}-
\Pi^h z_\delta\right)dx
\end{eqnarray*} 
will be taken into account. We then consider a  {\it unique} minimizer $\left( Q^h, f^h, g^h \right)$
of the {\it strictly convex} problem
\begin{eqnarray}\label{29-9-16ct1}
\min_{(Q,f,g) \in \mathcal{H}_{ad}}
{\mathcal{J}}^h_\delta(Q,f,g) + \rho\mathcal{R}(Q,f,g)
\end{eqnarray}
as a discrete regularized solution of the identification problem.  Note that, by using variational discretization concept introduced in \cite{hinze}, every solution of the minimization problem (\ref{29-9-16ct1}) is proved to automatically belong to finite dimensional spaces. Thus, a discretization of the admissible set $\mathcal{H}_{ad}$ can be avoided. Furthermore,  for simplicity of exposition we here restrict
ourselves to the case of one set of data $(z_\delta)_{\delta>0}$. In case with several sets of data $(z_{\delta_i})_{i=1}^{I}$ being available, we can replace the misfit term in the problem (\ref{29-9-16ct1}) by the term $\dfrac{1}{I}\sum_{i=1}^I{\mathcal{J}}^h_{\delta_i}(Q,f,g)$.

In \S \ref{Stability} we will show the convergence of these approximation solutions $\left( Q^h, f^h, g^h \right)$ to the identification $\left( Q^\dag,f^\dag,g^\dag\right)$ in the
$\mathbf{L}^2_{ \mbox{\tiny sym}}(\Omega) \times L^2(\Omega) \times L^2(\partial\Omega)$-norm as well as the convergence of corresponding approximation states $\big(\mathcal{U}^h_{Q^h, f^h, g^h}\big)$ to the exact $\Phi^\dag$ in the $H^1(\Omega)$-norm. Under the structural source condition --- but without the smallness requirement --- of the general convergence theory for non-linear, ill-posed problems (cf.\ \cite{Engl_Hanke_Neubauer,EnglKuNe}), we prove in \S \ref{Convergence rates} error bounds for these discrete approximations. For the numerical solution of the minimization problem (\ref{29-9-16ct1}) we in \S \ref{iterative} employ a gradient projection algorithm with Armijo steplength rule. Finally, a numerical implementation will be performed to illustrate  the theoretical findings.

The coefficient identification problem in PDEs arises from different contexts of applied sciences, e.g., from aquifer
analysis, geophysical prospecting and pollutant detection, and attracted great attention from many scientists in the last 30 years or so. For surveys on the subject one may consult in \cite{banks-kunisch1989,chavent2009,isakov,schuster,sun,Tarantola}. 
The problem of identifying the {\it scalar diffusion} coefficient has been extensively studied  for both theoretical research and numerical implementation, see e.g., \cite{ChanTai2003,ChanTai2004,Chavent_Kunisch2002,Chicone,Falk,hanke,Haoq,ito-kunisch,Kaltenbacher_Schoberl,keung-zou,knowles,kolo,Ric,wang_zou}. Some contributions for the case of the {\it simultaneous identification} can be found in \cite{Baku,hao_quyen3,hein,Know3} while some works treated the {\it diffusion matrix} case have been obtained in \cite{Deckelnick,Hinze-Tran,Hoffmann_Sprekels,Hsiao,Rannacher_Vexler}.

We conclude this introduction with the following mention. By using the H-convergent concept, the convergence analysis presented in \cite{Hinze-Tran} can not be applied directly to the problem of identifying {\it scalar diffusion} coefficients. There are two main difficulties for the scalar coefficient identification. First, the set
\begin{eqnarray*}
\mathcal{D} &:=& \left\{
qI_d \quad\big|\quad q\in L^\infty(\Omega) \quad\mbox{with}\quad \underline{q}\le q(x) \le\overline{q} \quad\mbox{a.e. in}\quad \Omega \quad\mbox{and}\quad I_d \quad\mbox{is the unit}\quad d\times d\mbox{-matrix}
\right\}
\end{eqnarray*}
is in general not a closed subset of $\mathcal{Q}_{ad}$ under the topology of the H-convergence (cf.\ \cite{tartar}), i.e. if the sequence $(q_nI_d)_n \subset \mathcal{D}$ is H-convergent to $Q\in\mathcal{Q}_{ad}$, then $Q$ is not necessarily proportional to $I_d$ in
dimension $d\ge 2$ or $Q\notin\mathcal{D}$. Second, the forward operator $\mathcal{U}$ is not weakly sequentially closed in $L^2$, i.e. if $(q_n,\mathcal{U}(q_n)) \rightharpoonup (q, \mathcal{Y})$ weakly in $L^2(\Omega)\times L^2(\Omega)$, it is not guaranteed that $\mathcal{Y}= \mathcal{U}(q)$ (see \cite{Deckelnick} and the references therein for counterexamples). To overcome these difficulties, a different analysis technique based on the convexity of the cost functional will be taken into counting. Due to the weak$^*$ closedness of the set $\mathcal{D}$ above in $\mathbf{L}^\infty_{ \mbox{\tiny sym}}(\Omega)$ (cf.\ Remark \ref{13-9-16ct1}), the convergence analysis performed in the present paper thus covers the scalar diffusion identification case. On the other hand, in \cite{Hinze-Tran} the source term and the boundary condition were assumed to be given. In the present situation they are variables which have to be found simultaneously together with the diffusion from observations.

Throughout the paper we write
$\int_\Omega \cdots$ instead of $\int_\Omega \cdots dx$ for the convenience of relevant notations. We use the
standard notion of Sobolev spaces $H^1(\Omega)$, $H^2(\Omega)$, $W^{k,p}(\Omega)$, etc from, e.g., \cite{Attouch}.

\section{Finite element discretization}\label{Finite element method}

\subsection{Preliminaries}\label{definition}

In product spaces $\mathbf{L}^2_{ \mbox{\tiny sym}}(\Omega) \times L^2(\Omega) \times L^2(\partial\Omega)$ and $\mathbf{L}^\infty_{ \mbox{\tiny sym}}(\Omega) \times L^2(\Omega)\times L^2(\partial\Omega)$ we use respectively the norm
\begin{eqnarray*}
\| (H,l,s) \|_{\mathbf{L}^2_{ \mbox{\tiny sym}}(\Omega) \times L^2(\Omega)\times L^2(\partial\Omega)} &=& \left(  \| H \|^2_{\mathbf{L}^2_{ \mbox{\tiny sym}}(\Omega)} + \| l \|^2_{L^2(\Omega)} + \|s\|^2_{L^2(\partial\Omega)}\right) ^{1/2} \mbox{~and}\\
\| (H,l,s) \|_{\mathbf{L}^\infty_{ \mbox{\tiny sym}}(\Omega) \times L^2(\Omega)\times L^2(\partial\Omega)} &=& \| H \|_{\mathbf{L}^\infty_{ \mbox{\tiny sym}}(\Omega)}+ \| l \|_{L^2(\Omega)} +\|s\|_{L^2(\partial\Omega)}.
\end{eqnarray*}

We note that the coefficient-to-solution operator 
\begin{eqnarray*}
\mathcal{U} : \quad\mathcal{H}_{ad} \subset
\mathbf{L}^\infty_{ \mbox{\tiny sym}}(\Omega)\times
L^2(\Omega)\times L^2(\partial\Omega) &\rightarrow& H^1_\diamond(\Omega)
\end{eqnarray*}
with
\begin{eqnarray*}
\Gamma := (Q,f,g)\in \mathcal{H}_{ad} &\rightarrow& \mathcal{U}(\Gamma) := \mathcal{U}_\Gamma
\end{eqnarray*}
is Fr\'echet differentiable on $\mathcal{H}_{ad}$. For each $\Gamma = (Q,
f,g)\in \mathcal{H}_{ad}$ the action of its Fr\'echet derivative in direction $\lambda := (H,l,s)\in \mathbf{L}^\infty_{ \mbox{\tiny sym}}(\Omega) \times L^2(\Omega)\times L^2(\partial\Omega)$ denoted by
$\xi_{\lambda}:=\mathcal{U}'_{\Gamma}(\lambda):=\mathcal{U}'(\Gamma)(\lambda)$ is
the unique weak solution in $H^1_\diamond(\Omega)$ to the equation
\begin{eqnarray}\label{22-9-16ct1}
\int_\Omega Q\nabla \xi_{\lambda} \cdot \nabla \varphi 
&=& -\int_\Omega H\nabla \mathcal{U}_{\Gamma} \cdot \nabla \varphi   +(l,\varphi) + \langle s, \gamma\varphi\rangle
\end{eqnarray}
for all $\varphi \in H^1_\diamond(\Omega)$. 

In $\mathcal{S}_d$ we introduce the convex subset
\begin{eqnarray*}
\mathcal{K} &:=& \big\{M \in \mathcal{S}_d \quad\big|\quad \underline{q}
\le M\xi\cdot\xi \le \overline{q} \quad\mbox{for all}\quad \xi \in R^d\big\}
\end{eqnarray*}
together with the orthogonal projection
$P_{\mathcal{K}} : \mathcal{S}_d \to \mathcal{K}$ that is characterised by
\begin{eqnarray*}
(A - P_{\mathcal{K}}(A)) \cdot (B - P_{\mathcal{K}}(A)) &\le& 0
\end{eqnarray*}
for all $A \in \mathcal{S}_d$ and $B \in \mathcal{K}$. Furthermore, let
$\xi := (\xi_1, \cdot\cdot\cdot, \xi_d)$ and $\eta :=
(\eta_1, \cdot\cdot\cdot, \eta_d)$ be two arbitrary vectors in $R^d$, we use
the notation
$$(\xi \otimes \eta)_{1\le i,j\le d} \in \mathcal{S}_d
\quad\mbox{with}\quad (\xi \otimes \eta)_{ij} := \frac{1}{2} (\xi_i \eta_j +
\xi_j \eta_i)  \quad\mbox{for all}\quad  i, j = 1, \cdots, d.$$

We close this subsection by the following note.
\begin{remark}\label{13-9-16ct1}
Let
\begin{eqnarray*}
\mathbf{D} &:=& \left\{
q\in L^\infty(\Omega) \quad\big|\quad \underline{q}\le q(x) \le\overline{q} \quad\mbox{a.e. in}\quad \Omega \right\}.
\end{eqnarray*}
Then $\mathbf{D}$ is a weakly$^*$ compact subset of $L^\infty(\Omega)$, i.e. for any sequence $(q_n)_n\subset\mathbf{D}$ a subsequence $(q_{n_m})_m$ and an element $\xi_\infty\in\mathbf{D}$ exist such that $(q_{n_m})_m$ is weakly$^*$ convergent in $L^\infty(\Omega)$ to $\xi_\infty$. In other words, for all $\theta_1\in L^1(\Omega)$ there holds the limit
\begin{eqnarray*}\lim_{m\to\infty} \int_\Omega
q_{n_m}\theta_1 &=& \int_\Omega \xi_\infty\theta_1.
\end{eqnarray*}
\end{remark}
We also remark that any $\Psi \in L^\infty(\Omega)$ can be considered as an element in ${L^\infty(\Omega)}^*$ by
\begin{eqnarray}
\left\<\Psi, \psi\right\>_{\big({L^\infty(\Omega)}^*,
L^\infty(\Omega)\big)} &:=& \int_\Omega \Psi \psi  \label{inf1}
\end{eqnarray}
for all $\psi \mbox{~in~} L^\infty(\Omega)$ and 
$\| \Psi \|_{{L^{\infty}(\Omega)}^*}\le |\Omega| \cdot \| \Psi \|_{L^\infty(\Omega)}$. Therefore, due to (\ref{inf1}), the assertion of Remark \ref{13-9-16ct1} is a direct consequence of the Banach-Alaoglu theorem.

\subsection{Discretization}

Let $\left(\mathcal{T}^h\right)_{0<h<1}$ be a family of regular and
quasi-uniform triangulations of the domain $\overline{\Omega}$ with the mesh size $h$ such that each vertex of the polygonal boundary $\partial\Omega$ is a node of $\mathcal{T}_h$.
For the definition of the discretization space of the state
functions let us denote
\begin{eqnarray}\label{27-10-15ct1}
\mathcal{V}^h_1 &=& \left\{\varphi^h\in C\left(\overline\Omega\right)\cap H^1_\diamond(\Omega)
\quad|\quad{\varphi^h}_{|T} \in \mathcal{P}_1(T) \quad\mbox{for all}\quad
T\in \mathcal{T}^h\right\}
\end{eqnarray}
with $\mathcal{P}_r$ consisting all polynomial functions of degree
at most $r$.
Similar to the continuous case, we have the following result.
\begin{lemma}
Let $(Q,f,g)$ be in $\mathcal{H}_{ad}$. Then the variational equation
\begin{eqnarray}
\int_\Omega Q\nabla \Phi^h \cdot \nabla \varphi^h &=&
(f,\varphi^h)+\langle g,\gamma\varphi^h\rangle\label{10/4:ct1}
\end{eqnarray}
for all $\varphi^h\in
\mathcal{V}^h_{1} $ admits a unique solution $\Phi^h \in \mathcal{V}^h_{1}$. Furthermore, the priori estimate
\begin{eqnarray}\label{18/5:ct1}
\|\Phi^h\|_{H^1(\Omega)} &\le& C_\mathcal{N}\left( \|f\|_{L^2(\Omega)} + \|g\|_{L^2(\partial\Omega)}\right)
\end{eqnarray}
is satisfied.
\end{lemma}

The map $\mathcal{U}^h : \mathcal{H}_{ad} \subset
\mathbf{L}^\infty_{ \mbox{\tiny sym}}(\Omega)\times
L^2(\Omega) \times
L^2(\partial\Omega) \rightarrow
\mathcal{V}^h_{1}$ from each $\Gamma := (Q,f,g) \in  \mathcal{H}_{ad}$
to the unique solution $\mathcal{U}^h_{\Gamma} := \Phi^h$ of (\ref{10/4:ct1})
is called the {\it discrete coefficient-to-solution
operator}.
This operator is also
Fr\'echet differentiable on the set $\mathcal{H}_{ad}$. For each $\Gamma = (Q,f,g)\in \mathcal{H}_{ad}$ and $\lambda := (H,l,s)\in \mathbf{L}^\infty_{ \mbox{\tiny sym}}(\Omega)\times L^2(\Omega) \times L^2(\partial\Omega)$ the Fr\'echet differential
$\xi^h_{\lambda}:={\mathcal{U}^{h}_{\Gamma}}'(\lambda)$ is an element of $\mathcal{V}^h_{1}$ and satisfies for all $\varphi^h$ in $\mathcal{V}^h_{1}$ the equation
\begin{eqnarray}\label{ct21}
\int_\Omega Q\nabla \xi^h_{\lambda} \cdot \nabla \varphi^h
&=& -\int_\Omega H\nabla \mathcal{U}^h_{\Gamma} \cdot \nabla \varphi^h  + (l,\varphi^h) + \langle s,\gamma \varphi^h\rangle.
\end{eqnarray}
Due to the standard theory of the finite element method for elliptic problems (cf.\ \cite{Brenner_Scott,ciarlet}), for any fixed $\Gamma = (Q,f,g) \in \mathcal{H}_{ad}$ there holds the limit
\begin{eqnarray}\label{24-10-15ct1}
\lim_{h \to 0} \left \| \mathcal{U}_{\Gamma} - \mathcal{U}^h_{\Gamma} \right\|_{H^1(\Omega)} \= 0.
\end{eqnarray}
Let 
\begin{eqnarray*}
\Pi^h: \quad L^1(\Omega) &\rightarrow& \left\{\varphi^h\in C\left(\overline\Omega\right)
\quad|\quad {\varphi^h}_{|T} \in \mathcal{P}_1(T) \quad\mbox{for all}\quad
T\in \mathcal{T}^h\right\} 
\end{eqnarray*}
be the Cl\'ement's mollification interpolation operator with properties
\begin{eqnarray}\label{23/10:ct2}
\lim_{h\to 0} \left\| \phi - \Pi^h \phi
\right\|_{H^k(\Omega)} \= 0 \quad\mbox{for all}\quad k \in \{0,
1\}
\end{eqnarray}
and
\begin{eqnarray}\label{23/5:ct1}
\left\| \phi - \Pi^h \phi \right\|_{H^k(\Omega)} &\le&
Ch^{l-k} \| \phi\|_{H^l(\Omega)}
\end{eqnarray}
for $0 \le k \le l \le 2$, where $C$ is independent of $h$ and $\phi$ (cf.\ \cite{Clement,Bernardi1,Bernardi2,scott_zhang}). Then, using the discrete operator $\mathcal{U}^h$ and the interpolation operator $\Pi^h$, we can now introduce
the discrete cost functional
\begin{eqnarray}\label{29/6:ct9}
\mathcal{J}^h_\delta (Q,f,g) &:=& \int_\Omega Q\nabla \left(\mathcal{U}^h_{Q,f,g}-
\Pi^h z_\delta\right) \cdot  \nabla \left(\mathcal{U}^h_{Q,f,g}-
\Pi^h z_\delta\right),
\end{eqnarray}
where $(Q,f,g) \in \mathcal{H}_{ad}$. 

\begin{lemma}\label{finite-dim}
Assume that the sequence $\left( \Gamma_n\right)_n := \left(Q_n, f_n, g_n\right)_n\subset \mathcal{H}_{ad}$
weakly converges to $\Gamma := (Q,f,g)$ in $\mathbf{L}^2_{ \mbox{\tiny sym}}(\Omega) \times
L^2(\Omega)\times L^2(\partial\Omega)$. Then for any fixed $h>0$ the sequence $\left(\mathcal{U}^h_{\Gamma_n}\right)_n \subset \mathcal{V}^h_{1}$ converges to $\mathcal{U}^h_{\Gamma}$ in the $H^1(\Omega)$-norm.
\end{lemma}

\begin{proof}
Due to Remark \ref{13-9-16ct1}, $(Q_n)_n$ has a subsequence denoted by the same symbol which is weakly$^*$ convergent in $\mathbf{L}^\infty_{ \mbox{\tiny sym}}(\Omega)$ to $Q$. Furthermore, by (\ref{18/5:ct1}),
the corresponding state sequence $\left(\mathcal{U}^h_{\Gamma_n}\right)_n$ is bounded in the {\it finite} dimensional space $\mathcal{V}^h_{1}$. A subsequence which is not relabelled and an element $\Theta^h \in \mathcal{V}^h_{1}$ then exist such that $\left(\mathcal{U}^h_{\Gamma_n}\right)_n$ converges to $\Theta^h$ in the $H^1(\Omega)$-norm. It follows from the equation (\ref{10/4:ct1}) that
\begin{eqnarray}\label{14-9-16ct1}
\int_\Omega Q_n \nabla \big(  \mathcal{U}^h_{\Gamma_n} -  \mathcal{U}^h_{\Gamma}\big) \cdot \nabla \varphi^h &=& \int_\Omega \left( Q - Q_n \right) \nabla \mathcal{U}^h_{\Gamma} \cdot \nabla \varphi^h  + \left( f_n - f,\varphi^h\right)  +\left\langle g_n -g,\gamma\varphi^h \right\rangle
\end{eqnarray}
for all $\varphi^h \in \mathcal{V}^h_{1}$. Taking $\varphi^h = \mathcal{U}^h_{\Gamma_n} -  \mathcal{U}^h_{\Gamma}$, by (\ref{coercivity}), we obtain that
\begin{eqnarray}\label{30-9-16ct3}
\frac{C_\Omega \underline{q}}{1+ C_\Omega} \left\|  \mathcal{U}^h_{\Gamma_n} -  \mathcal{U}^h_{\Gamma} \right \|^2_{H^1(\Omega)}
&\le& \int_\Omega \left( Q - Q_n \right) \nabla \mathcal{U}^h_{\Gamma} \cdot \nabla \left( \mathcal{U}^h_{\Gamma_n} -  \Theta^h + \Theta^h - \mathcal{U}^h_{\Gamma}\right)   \\
&~\quad& + \left( f_n - f, \mathcal{U}^h_{\Gamma_n} -  \Theta^h + \Theta^h -  \mathcal{U}^h_{\Gamma}\right)  + \left\langle g_n - g, \gamma \left( \mathcal{U}^h_{\Gamma_n} -  \Theta^h + \Theta^h -  \mathcal{U}^h_{\Gamma}\right) \right\rangle \nonumber\\
&\le& C \left\|  \mathcal{U}^h_{\Gamma_n} -  \Theta^h \right \|_{H^1(\Omega)} + \int_\Omega \left( Q - Q_n \right) \nabla \mathcal{U}^h_{\Gamma} \cdot \nabla \left( \Theta^h - \mathcal{U}^h_{\Gamma}\right)  \nonumber\\
&~\quad& + \left( f_n - f, \Theta^h -  \mathcal{U}^h_{\Gamma}\right)  + \left\langle g_n - g, \gamma \left(  \Theta^h -  \mathcal{U}^h_{\Gamma}\right) \right\rangle.\nonumber
\end{eqnarray}
Since $Q_n\rightharpoonup Q$ weakly$^*$ in $\mathbf{L}^\infty_{ \mbox{\tiny sym}}(\Omega)$, 
we get
$\limn \int_\Omega \left( Q - Q_n \right) \nabla \mathcal{U}^h_{\Gamma} \cdot \nabla \left( \Theta^h - \mathcal{U}^h_{\Gamma}\right) =0.$
Sending $n$ to $\infty$, we thus obtain from the last inequality that $\limn \left\|  \mathcal{U}^h_{\Gamma_n} -  \mathcal{U}^h_{\Gamma} \right \|_{H^1(\Omega)} = 0$,
which finishes the proof.
\end{proof}

We now state the following useful result on the convexity of the cost functional.

\begin{lemma} \label{J-dis-convex}
$\mathcal{J}^h_\delta$
is convex and continuous on
$\mathcal{H}_{ad}$ with respect to the $\mathbf{L}^2_{ \mbox{\tiny sym}}(\Omega) \times
L^2(\Omega)\times L^2(\partial\Omega)$-norm.
\end{lemma}

\begin{proof}
The continuity of $\mathcal{J}^h_\delta$ follows directly from Lemma \ref{finite-dim}. We show that $\mathcal{J}^h_\delta$ is convex.

Let $\Gamma := (Q,f,g) \in \mathcal{H}_{ad}$ and $\lambda := (H,l,s)\in
\mathbf{L}^{\infty}_{ \mbox{\tiny sym}}(\Omega)\times L^2(\Omega) \times L^2(\partial\Omega)$. We have that
\begin{eqnarray*}
{\mathcal{U}^h_{\Gamma}}'(\lambda) \= \frac{\partial \mathcal{U}^h_{\Gamma}}{\partial
Q} H +
\frac{\partial \mathcal{U}^h_{\Gamma}}{\partial f} l + \frac{\partial \mathcal{U}^h_{\Gamma}}{\partial g} s
\quad\mbox{and}\quad
{\mathcal{J}^h_\delta}'(\Gamma)(\lambda) \quad=\quad \frac{\partial
\mathcal{J}^h_\delta(\Gamma)}{\partial Q} H   + \frac{\partial \mathcal{J}^h_\delta(\Gamma)}{\partial f} l + \frac{\partial \mathcal{J}^h_\delta(\Gamma)}{\partial g} s.
\end{eqnarray*}
We compute for each term in the right hand side of the last equation. First we get
\begin{eqnarray*}
\frac{\partial \mathcal{J}^h_\delta(\Gamma)}{\partial Q} H
&=& \int_\Omega  H\nabla \left(\mathcal{U}^h_{\Gamma} - \Pi^hz_\delta \right) \cdot
\nabla \left(\mathcal{U}^h_{\Gamma} - \Pi^hz_\delta \right) + 2\int_\Omega Q \nabla \left(\frac{\partial \mathcal{U}^h_{\Gamma}}{\partial Q}H\right) \cdot \nabla \left(\mathcal{U}^h_{\Gamma} - \Pi^hz_\delta \right).
\end{eqnarray*}
For the second term we have
\begin{eqnarray*}
\frac{\partial \mathcal{J}^h_\delta(\Gamma)}{\partial f} l
\= 2 \int_\Omega  Q \nabla \left(\frac{\partial \mathcal{U}^h_{\Gamma}}{\partial
f}l \right)  \cdot \nabla \left(\mathcal{U}^h_{\Gamma} - \Pi^hz_\delta \right).
\end{eqnarray*}
Finally, we have for the third term
\begin{eqnarray*}
\frac{\partial \mathcal{J}^h_\delta(\Gamma)}{\partial g} s
\= 2 \int_\Omega  Q \nabla \left(\frac{\partial \mathcal{U}^h_{\Gamma}}{\partial g}s \right)  \cdot \nabla \left(\mathcal{U}^h_{\Gamma} - \Pi^hz_\delta \right).
\end{eqnarray*}
Therefore,
\begin{eqnarray*}
{\mathcal{J}^h_\delta}'(\Gamma)(\lambda) &=& 2\int_\Omega Q \nabla \left(\frac{\partial \mathcal{U}^h_{\Gamma}}{\partial Q}H + \frac{\partial \mathcal{U}^h_{\Gamma}}{\partial f}l + \frac{\partial \mathcal{U}^h_{\Gamma}}{\partial g}s\right) \cdot \nabla \left(\mathcal{U}^h_{\Gamma} - \Pi^hz_\delta \right)  + \int_\Omega H\nabla
(\mathcal{U}^h_{\Gamma} - \Pi^hz_\delta ) \cdot \nabla (\mathcal{U}^h_{\Gamma} -
\Pi^hz_\delta ) \\
&=& 2\int_\Omega Q \nabla {\mathcal{U}^h_{\Gamma}}'(\lambda) \cdot \nabla \left(\mathcal{U}^h_{\Gamma} - \Pi^hz_\delta \right) 
+ \int_\Omega H\nabla \left(\mathcal{U}^h_{\Gamma} - \Pi^hz_\delta \right) \cdot \nabla \left(\mathcal{U}^h_{\Gamma} - \Pi^hz_\delta \right)\\
&=& 2\int_\Omega Q \nabla {\mathcal{U}^h_{\Gamma}}'(\lambda) \cdot \nabla \left(\mathcal{U}^h_{\Gamma} - \bar{\Pi}^hz_\delta \right) 
+ \int_\Omega H\nabla \left(\mathcal{U}^h_{\Gamma} - \bar{\Pi}^hz_\delta \right) \cdot \nabla \left(\mathcal{U}^h_{\Gamma} - \bar{\Pi}^hz_\delta \right),
\end{eqnarray*}
where
\begin{eqnarray}\label{21-9-16ct4}
\bar{\Pi}^hz_\delta &:=& \Pi^hz_\delta - |\Omega|^{-1} \left\langle 1, \gamma\Pi^hz_\delta \right\rangle  \in\mathcal{V}^h_1 \quad\mbox{with}\quad \nabla \bar{\Pi}^hz_\delta \quad=\quad \nabla\Pi^hz_\delta.
\end{eqnarray}
By (\ref{ct21}), we infer that
\begin{eqnarray}\label{23-10-15ct1}
{\mathcal{J}^h_\delta}'(\Gamma)(\lambda) 
&=& - 2\int_\Omega H \nabla \mathcal{U}^h_{\Gamma} \cdot \nabla
\left(\mathcal{U}^h_{\Gamma}-\bar{\Pi}^hz_\delta\right) + 2\left(l,\mathcal{U}^h_{\Gamma}-\bar{\Pi}^hz_\delta\right) +2 \left\langle s, \gamma\left( \mathcal{U}^h_{\Gamma}-\bar{\Pi}^hz_\delta\right)  \right\rangle \nonumber\\
&~\quad& + \int_\Omega H\nabla \left(\mathcal{U}^h_{\Gamma} - \bar{\Pi}^hz_\delta \right) \cdot \nabla \left(\mathcal{U}^h_{\Gamma} - \bar{\Pi}^hz_\delta \right)\\
&=& -\int_\Omega H\nabla \mathcal{U}^h_{\Gamma} \cdot \nabla
\mathcal{U}^h_{\Gamma} + \int_\Omega H\nabla \bar{\Pi}^hz_\delta \cdot \nabla \bar{\Pi}^hz_\delta + 2\left(l,\mathcal{U}^h_{\Gamma}-\bar{\Pi}^hz_\delta\right) +2 \left\langle s, \gamma\left( \mathcal{U}^h_{\Gamma}-\bar{\Pi}^hz_\delta\right)  \right\rangle. \nonumber
\end{eqnarray}
Therefore, by (\ref{ct21}) again, we arrive at
\begin{eqnarray*}
{\mathcal{J}^h_\delta}''(\Gamma)\left(\lambda, \lambda\right) &=&
-2\int_\Omega H \nabla \mathcal{U}^h_{\Gamma} \cdot \nabla {\mathcal{U}^h_{\Gamma}}'(\lambda) + 2 \left(  l, {\mathcal{U}^h_{\Gamma}}'(\lambda)\right) + 2 \left\langle s, \gamma{\mathcal{U}^h_{\Gamma}}'(\lambda) \right\rangle\\
&=& 2\int_\Omega Q\nabla {\mathcal{U}^h_{\Gamma}}'(\lambda) \cdot \nabla
{\mathcal{U}^h_{\Gamma}}'(\lambda)
\quad\ge\quad 2\frac{C_\Omega \underline{q}}{1+ C_\Omega} \left\| {\mathcal{U}^h_{\Gamma}}'(\lambda) \right\|^2_{H^1(\Omega)}  \quad\ge\quad 0,
\end{eqnarray*}
by (\ref{coercivity}), which completes the proof.
\end{proof}

Now we are in position to prove the main result of this section.

\begin{theorem} \label{solution2}
The strictly convex minimization problem
$$
\min_{(Q,f,g)\in \mathcal{H}_{ad}} \Upsilon^{\rho,h}_\delta (Q,f,g) \quad:=\quad \mathcal{J}^h_\delta(Q,f,g) +\rho\mathcal{R}(Q,f,g) \eqno \left(
\mathcal{P}^{\rho,h}_\delta \right) 
$$
attains a unique minimizer. Furthermore, an element $\Gamma:= (Q,f,g) \in \mathcal{H}_{ad}$ is the unique minimizer to $\left(
\mathcal{P}^{\rho,h}_\delta \right)$ if and only if the system
\begin{eqnarray}
Q(x) &=& P_{\mathcal{K}} \left( \frac{1}{2\rho}\Big(\nabla \mathcal{U}^h_\Gamma(x) \otimes \nabla \mathcal{U}^h_\Gamma(x) -\nabla\bar{\Pi}^h z_\delta(x) \otimes \nabla \bar{\Pi}^h z_\delta(x)\Big)\right), \label{21-9-16ct1}\\
f(x) &=& \frac{1}{\rho} \Big(\bar{\Pi}^h z_\delta (x) - \mathcal{U}^h_\Gamma(x)\Big),\label{21-9-16ct1*}\\
g(x) &=& \frac{1}{\rho} \gamma \Big(\bar{\Pi}^h z_\delta (x) - \mathcal{U}^h_\Gamma(x)\Big)\label{21-9-16ct1**}
\end{eqnarray}
holds for a.e. in $\Omega$, where $\bar{\Pi}^h$ was generated from $\Pi^h$ according to (\ref{21-9-16ct4}).
\end{theorem}

\begin{proof}
Let $(\Gamma_n)_n := (Q_n,f_n,g_n)_n \subset\mathcal{H}_{ad}$ be a minimizing sequence of $\left(
\mathcal{P}^{\rho,h}_\delta \right)$, i.e.
\begin{eqnarray*}
\limn \Upsilon^{\rho,h}_\delta(\Gamma_n) \= \inf_{(Q,f,g)\in \mathcal{H}_{ad}} \Upsilon^{\rho,h}_\delta (Q,f,g).
\end{eqnarray*}
The sequence $(\Gamma_n)_n$ is thus bounded in the $\mathbf{L}^2_{ \mbox{\tiny sym}}(\Omega) \times
L^2(\Omega)\times L^2(\partial\Omega)$-norm. A subsequence not relabelled and an element $\Gamma:= (Q,f,g)\in \mathbf{L}^2_{ \mbox{\tiny sym}}(\Omega) \times L^2(\Omega)\times L^2(\partial\Omega)$ exist such that $\Gamma_n \rightharpoonup \Gamma$ weakly in $\mathbf{L}^2_{ \mbox{\tiny sym}}(\Omega) \times L^2(\Omega)\times L^2(\partial\Omega)$. On the other hand, since $\mathcal{H}_{ad}$ is a convex, closed subset of $\mathbf{L}^2_{ \mbox{\tiny sym}}(\Omega) \times L^2(\Omega)\times L^2(\partial\Omega)$, so is weakly closed, it follows that $\Gamma\in\mathcal{H}_{ad}$.
By Lemma \ref{J-dis-convex}, $\mathcal{J}^h_\delta$ and $\mathcal{R}$ are both weakly lower semi-continuous on $\mathcal{H}_{ad}$ which yields that
\begin{eqnarray*}\mathcal{J}^h_\delta(\Gamma) &\le& \liminfn \mathcal{J}^h_\delta(\Gamma_n) \quad\mbox{and}\quad \mathcal{R}(\Gamma) \quad\le\quad\liminfn \mathcal{R}(\Gamma_n).
\end{eqnarray*}
We therefore have that
\begin{eqnarray*}
\mathcal{J}^h_\delta(\Gamma)+ \mathcal{R}(\Gamma) &\le& \liminfn \mathcal{J}^h_\delta(\Gamma_n)+ \liminfn \mathcal{R}(\Gamma_n) \quad\le\quad \liminfn \left( \mathcal{J}^h_\delta(\Gamma_n)+ \mathcal{R}(\Gamma_n)\right) \\
&=& \limn \Upsilon^{\rho,h}_\delta(\Gamma_n) \quad=\quad \inf_{(Q,f,g)\in \mathcal{H}_{ad}} \Upsilon^{\rho,h}_\delta (Q,f,g),
\end{eqnarray*}
and $\Gamma$ is then a minimizer to $\left(
\mathcal{P}^{\rho,h}_\delta \right)$. Since $\Upsilon^{\rho,h}_\delta$ is strictly convex, this minimizer is unique. Next, an element $\Gamma:= (Q,f,g) \in \mathcal{H}_{ad}$ is the minimizer to $\left(
\mathcal{P}^{\rho,h}_\delta \right)$ if and only if 
${\Upsilon^{\rho,h}_\delta}'(\Gamma)(\overline{\Gamma} - \Gamma) \ge 0$
for all $\overline{\Gamma} =(H,l,s)\in \mathcal{H}_{ad}$. Then, in view of (\ref{23-10-15ct1}), we get that
\begin{eqnarray*}
0& \le& \int_\Omega (H-Q)\nabla \bar{\Pi}^hz_\delta \cdot \nabla \bar{\Pi}^hz_\delta -\int_\Omega (H-Q)\nabla \mathcal{U}^h_{\Gamma} \cdot \nabla
\mathcal{U}^h_{\Gamma} +2\rho (H-Q,Q)\nonumber\\
&~\quad& + 2\left(l-f,\mathcal{U}^h_{\Gamma}-\bar{\Pi}^hz_\delta\right) +2\rho(l-f,f)+2 \left\langle s-g, \gamma\left( \mathcal{U}^h_{\Gamma}-\bar{\Pi}^hz_\delta\right)  \right\rangle +2\rho \left\langle s-g,g\right\rangle \nonumber\\
&=& \int_\Omega (H-Q) \cdot \left( \nabla \bar{\Pi}^hz_\delta \otimes \nabla \bar{\Pi}^hz_\delta - \nabla \mathcal{U}^h_{\Gamma} \otimes \nabla
\mathcal{U}^h_{\Gamma} +2\rho Q \right) \nonumber\\
&~\quad& + 2\left(l-f,\mathcal{U}^h_{\Gamma}-\bar{\Pi}^hz_\delta +\rho f\right) + 2 \left\langle s-g, \gamma\left( \mathcal{U}^h_{\Gamma}-\bar{\Pi}^hz_\delta\right) +\rho g \right\rangle
\end{eqnarray*}
for all $\overline{\Gamma} =(H,l,s)\in \mathcal{H}_{ad}$. Taking $\overline{\Gamma}_1 = (H,f,g)$, $\overline{\Gamma}_2 = (Q,l,g)$ and $\overline{\Gamma}_3 = (Q,f,s)$ into the above inequality we obtain the system (\ref{21-9-16ct1})--(\ref{21-9-16ct1**}). The proof is completed.
\end{proof}

\begin{remark}\label{22-9-16ct2}
We denote by
\begin{eqnarray*}
\mathcal{V}^h_0 &:=& \left\{\varphi^h\in L^2(\Omega)
\quad\big|\quad {\varphi^h}_{|T} \quad=\quad \mbox{const for all triangulations}\quad
T\in \mathcal{T}^h\right\},\\
\mathcal{E}^h_1 &:=& \left\{\varphi^h\in C(\partial\Omega)
\quad\big|\quad {\varphi^h}_{|e} \in \mathcal{P}_1 \quad\mbox{for all boundary edges}
\enskip e \enskip \mbox{of} \quad\mathcal{T}^h\right\}.
\end{eqnarray*}
Since $\mathcal{U}^h_{\Gamma}\in \mathcal{V}^h_1$ and $\bar{\Pi}^hz_\delta\in \mathcal{V}^h_1$, the system (\ref{21-9-16ct1})--(\ref{21-9-16ct1**}) shows that every solution of 
$\big(
\mathcal{P}^{\rho,h}_\delta \big)$ automatically belongs  to the finite dimensional space ${\mathcal{V}^h_0}^{d\times d} \times \mathcal{V}^h_1 \times \mathcal{E}^h_1$. 
\end{remark}

\section{Convergence}\label{Stability}

For abbreviation in what follows we denote by $C$ a generic positive constant independent of the mesh size $h$, the noise level $\delta$ and the regularization parameter $\rho$. By (\ref{23/10:ct2}) and (\ref{23/5:ct1}), we can introduce for each $\Phi \in H^1(\Omega)$
\begin{eqnarray*}
\chi^h_\Phi &:=& \left\| \Phi - \Pi^h\Phi \right\|_{H^1(\Omega)} \quad\mbox{which satisfies}\quad \lim_{h \to 0} \chi^h_\Phi \quad=\quad 0 \quad\mbox{and}\quad 0 \quad\le\quad \chi^h_\Phi \quad\le\quad Ch
\end{eqnarray*}
in case $\Phi \in H^2(\Omega)$. Likewise, by (\ref{24-10-15ct1}), for all $\Gamma\in\mathcal{H}_{ad}$
\begin{eqnarray*}
\beta^h_{\mathcal{U}_\Gamma} &:=& \left \| \mathcal{U}_{\Gamma} - \mathcal{U}^h_{\Gamma} \right\|_{H^1(\Omega)} \to 0 \quad\mbox{as}\quad h \to 0 \quad\mbox{and}\quad 0 \quad\le\quad \beta^h_{\mathcal{U}_\Gamma} \quad\le\quad Ch \quad\mbox{as}\quad \mathcal{U}_\Gamma\in H^2(\Omega).
\end{eqnarray*}
Furthermore, by (\ref{23/5:ct1}), we get
\begin{eqnarray}\label{15-9-16ct3}
\|\Pi^h\|_{\mathcal{L}\left(L^2(\Omega),L^2(\Omega)\right)} &\le& C \quad\mbox{and}\quad \|\Pi^h\|_{\mathcal{L}\left(H^1(\Omega),H^1(\Omega)\right)} \quad\le\quad C.
\end{eqnarray}
Thus, it follows from the inverse inequality (cf.\ \cite{Brenner_Scott,ciarlet}): 
\begin{eqnarray*}
\|\varphi^h\|_{H^1(\Omega)} &\le& Ch^{-1} \|\varphi^h\|_{L^2(\Omega)} \quad\mbox{for all}\quad \varphi^h \in \left\{\varphi^h\in C\left(\overline\Omega\right)
\quad\big|\quad {\varphi^h}_{|T} \in \mathcal{P}_1(T) \quad\mbox{for all}\quad
T\in \mathcal{T}^h\right\}
\end{eqnarray*}
that
\begin{eqnarray}\label{21-9-16ct3}
\|\Phi^\dag -\Pi^hz_\delta\|_{H^1(\Omega)} 
&\le& \|\Pi^h \left( \Phi^\dag -z_\delta\right) \|_{H^1(\Omega)} + \|\Phi^\dag -\Pi^h\Phi^\dag\|_{H^1(\Omega)} \quad\le\quad Ch^{-1}\|\Pi^h \left( \Phi^\dag -z_\delta\right) \|_{L^2(\Omega)}+\chi^h_{\Phi^\dag} \nonumber\\
&\le& Ch^{-1} \|\Pi^h\|_{\mathcal{L}\left(L^2(\Omega),L^2(\Omega)\right)}\|\Phi^\dag -z_\delta \|_{L^2(\Omega)}+\chi^h_{\Phi^\dag} \quad\le\quad Ch^{-1}\delta+\chi^h_{\Phi^\dag}.
\end{eqnarray}
The following result shows the convergence of
finite element approximations to the unique minimum norm solution $\Gamma^\dag := \left(Q^\dag, f^\dag, g^\dag\right)$ of the identification problem, which is defined by (\ref{8-1-16ct2}). 

\begin{theorem}\label{Stability2}
Let $\left( h_n\right) _n$ be a sequence with $\limn h_n = 0$ and $\left( \delta_n\right) _n$ and $\left( \rho_n\right) _n$ are any positive sequences such that
\begin{eqnarray*}\rho_n \rightarrow 0, \quad\frac{\delta_n}{h_n\sqrt{\rho_n}}
\rightarrow 0, \quad\frac{\beta^{h_n}_{\mathcal{U}_{\Gamma^\dag}}}{\sqrt{\rho_n}} \rightarrow 0 \quad\mbox{and}\quad \frac{\chi^{h_n}_{\Phi^\dag}}{\sqrt{\rho_n}} \rightarrow 0 \quad\mbox{as}\quad n \to \infty.
\end{eqnarray*}
Assume that $\left( z_{\delta_n}\right)_n \subset L^2(\Omega)$ is a sequence satisfying $\left\| z_{\delta_n} - \Phi^\dag\right\|_{L^2(\Omega)} \le \delta_n$ and $\Gamma_n := \left( Q_n,f_n,g_n \right)$ 
is the unique minimizer of the problem $\big( \mathcal{P}^{\rho_n, h_n}_{\delta_n} \big)$ for each $n\in N$. Then the sequence $\left( \Gamma_n\right)_n$
converges to $\Gamma^\dag$ in the
$\mathbf{L}^2_{ \mbox{\tiny sym}}(\Omega) \times L^2(\Omega) \times L^2(\partial\Omega)$-norm as $n \to \infty$. Furthermore, the corresponding discrete state sequence $\big(\mathcal{U}^{h_n}_{\Gamma_n}\big)_n$ also converges to $\Phi^\dag$ in the $H^1(\Omega)$-norm.
\end{theorem}

\begin{remark}\label{15-9-16ct1}
In case $\Phi^\dag = \mathcal{U}_{\Gamma^\dag} \in H^2(\Omega)$ we have $0\le \beta^{h_n}_{\mathcal{U}_{\Gamma^\dag}}, \chi^{h_n}_{\Phi^\dag} \le Ch_n$. Therefore, the convergence of Theorem \ref{Stability2} is obtained if $\delta_n \sim h_n^2$ and the sequence $(\rho_n)_n$ is chosen such that
$$\rho_n \rightarrow 0 \quad\mbox{and}\quad \frac{h_n}{\sqrt{\rho_n}}
\rightarrow 0 \quad\mbox{as}\quad n \to \infty.$$
\end{remark}

To prove Theorem \ref{Stability2}, we need the following auxiliary estimate.

\begin{lemma}\label{8-1-16ct1}
There holds the estimate
\begin{eqnarray}\label{8-1-16ct4*}
\mathcal{J}_{\delta}^{h} (\Gamma^\dag) &\le& C \left( h^{-2}\delta^2+ \big(\chi^h_{\Phi^\dag}\big)^2 + \big(\beta^h_{\mathcal{U}_{\Gamma^\dag}}\big)^2\right).
\end{eqnarray}
\end{lemma}

\begin{proof}
We have with $\Phi^\dag = \mathcal{U}_{\Gamma^\dag}$ and (\ref{21-9-16ct3}) that
\begin{eqnarray*}
\mathcal{J}_{\delta}^{h} (\Gamma^\dag) 
&=& \int_\Omega Q^\dag \nabla \big(\mathcal{U}^h_{\Gamma^\dag} -\Pi^hz_\delta\big) \cdot \nabla \big(\mathcal{U}^h_{\Gamma^\dag} -\Pi^hz_\delta\big)
\quad\le\quad \overline{q} \left\|  \mathcal{U}^h_{\Gamma^\dag} -   \Pi^hz_\delta \right \|^2_{H^1(\Omega)}\nonumber\\
&=& \overline{q} \left\|  \mathcal{U}^h_{\Gamma^\dag} - \mathcal{U}_{\Gamma^\dag} + \Phi^\dag -   \Pi^hz_\delta \right \|^2_{H^1(\Omega)} \quad\le\quad C \left( \left\|  \mathcal{U}^h_{\Gamma^\dag} - \mathcal{U}_{\Gamma^\dag} \right \|^2_{H^1(\Omega)} + \left\|  \Phi^\dag -   \Pi^hz_\delta \right \|^2_{H^1(\Omega)}\right) \\
&\le& C \left( h^{-2}\delta^2+ \big(\chi^h_{\Phi^\dag}\big)^2 + \big(\beta^h_{\mathcal{U}_{\Gamma^\dag}}\big)^2\right),
\end{eqnarray*}
which finishes the proof.
\end{proof}

\begin{proof}[Proof of Theorem \ref{Stability2}]
By the optimality of $\Gamma_n$ and Lemma \ref{8-1-16ct1}, we have that
\begin{eqnarray*}
\mathcal{J}_{\delta_n}^{h_n} \left(\Gamma_n\right) + \rho_n \mathcal{R}
\left(\Gamma_n\right)
&\le& \mathcal{J}_{\delta_n}^{h_n} \left( \Gamma^\dag \right) + \rho_n
\mathcal{R} \left( \Gamma^\dag\right) \nonumber\\
&\le& C \left( h_n^{-2}\delta_n^2+ \big(\chi^{h_n}_{\Phi^\dag}\big)^2 + \big(\beta^{h_n}_{\mathcal{U}_{\Gamma^\dag}}\big)^2\right)  + \rho_n \mathcal{R} \left( \Gamma^\dag\right)
\end{eqnarray*}
which yields
\begin{eqnarray}\label{odinh22}
\limn \mathcal{J}_{\delta_n}^{h_n} \left(\Gamma_n\right) \= 0
\end{eqnarray}
and
\begin{eqnarray}
\limsupn \mathcal{R}\left(\Gamma_n\right)
&\le&  \mathcal{R}\left(\Gamma^\dag\right).\label{odinh22*}
\end{eqnarray}
A subsequence of the sequence $\left(\Gamma_n\right)_n$ denoted by the same symbol and an element
$\Gamma_0 := (Q_0, f_0, g_0) \in \mathcal{H}_{ad} $
then exist such that
\begin{eqnarray*}
Q_n &\rightharpoonup& Q_0 \quad\mbox{weakly* in}\quad \mathbf{L}^\infty_{ \mbox{\tiny sym}}(\Omega),\\
f_n &\rightharpoonup& f_0 \quad\mbox{weakly in}\quad L^2(\Omega),\\
g_n &\rightharpoonup& g_0 \quad\mbox{weakly in}\quad L^2(\partial\Omega).
\end{eqnarray*}
We will show that $\left(\Gamma_n\right)_n$
converges to $\Gamma_0$ in the $\mathbf{L}^2_{ \mbox{\tiny sym}}(\Omega) \times L^2(\Omega) \times L^2(\partial\Omega)$-norm and $\Gamma_0 = \Gamma^\dag$. We have from (\ref{21-9-16ct3}) that
\begin{eqnarray}\label{2-12-15ct4}
\limn \left\| \Pi^{h_n} z_{\delta_n} - \mathcal{U}_{\Gamma^\dag} \right\|_{H^1(\Omega)} &\le& \limn \left( Ch^{-1}_n\delta_n + \chi^{h_n}_{\Phi^\dag}\right) \quad=\quad 0.
\end{eqnarray}
Combining this with $\limn \|\mathcal{U}_{\Gamma_0}-\mathcal{U}^{h_n}_{\Gamma_0}\|_{H^1(\Omega)} =0$ from (\ref{24-10-15ct1}), we arrive at
\begin{eqnarray*}
\limn \mathcal{J}^{h_n}_{\delta_n}
\left(\Gamma_0\right) 
&=& \limn \int_\Omega Q_0 \nabla \left( \mathcal{U}^{h_n}_{\Gamma_0} -\Pi^{h_n} z_{\delta_n}\right) \cdot \nabla \left( \mathcal{U}^{h_n}_{\Gamma_0} -\Pi^{h_n} z_{\delta_n}\right) \\
&=& \int_\Omega Q_0 \nabla \left( \mathcal{U}_{\Gamma_0} -\mathcal{U}_{\Gamma^\dag}\right) \cdot \nabla \left( \mathcal{U}_{\Gamma_0} -\mathcal{U}_{\Gamma^\dag}\right).
\end{eqnarray*}
Now for each fixed $n$ we consider an arbitrary subsequence $(\Gamma_{n_m})_m$ of $(\Gamma_n)_n$. By the weakly l.s.c. property of the functional $\mathcal{J}^{h_n}_{\delta_n}$ (cf.\ Lemma \ref{J-dis-convex}), we obtain that
\begin{eqnarray*}
\mathcal{J}^{h_n}_{\delta_n}
\left(\Gamma_0\right) &\le& \liminf_{m\to\infty} \mathcal{J}^{h_n}_{\delta_n} (\Gamma_{n_m}).
\end{eqnarray*}
Again, using the convexity of $\mathcal{J}^{h_n}_{\delta_n}$, we get that
\begin{eqnarray*}
\mathcal{J}^{h_n}_{\delta_n}
\left(\Gamma_n\right) &\ge& \mathcal{J}^{h_n}_{\delta_n} (\Gamma_{n_m}) +  {\mathcal{J}^{h_n}_{\delta_n}}' \left(\Gamma_{n_m}\right) \left( \Gamma_n - \Gamma_{n_m}\right).
\end{eqnarray*}
By (\ref{coercivity}), we thus arrive at
\begin{eqnarray*}
C \left\|\mathcal{U}_{\Gamma_0} -
\mathcal{U}_{\Gamma^\dag}\right\|^2_{H^1(\Omega)}
&\le& \int_\Omega Q_0 \nabla \left( \mathcal{U}_{\Gamma_0} -\mathcal{U}_{\Gamma^\dag}\right) \cdot \nabla \left( \mathcal{U}_{\Gamma_0} -\mathcal{U}_{\Gamma^\dag}\right) \\
&=&  \limn \mathcal{J}^{h_n}_{\delta_n}
\left(\Gamma_0\right) \quad\le\quad \limn \left( \liminf_{m\to\infty} \mathcal{J}^{h_n}_{\delta_n}
\left(\Gamma_{n_m}\right) \right) \\
&\le& \limn \liminf_{m\to\infty} \left( \mathcal{J}^{h_n}_{\delta_n}
\left(\Gamma_n\right) + {\mathcal{J}^{h_n}_{\delta_n}}' \left(\Gamma_{n_m}\right) \left( \Gamma_{n_m} - \Gamma_n\right) \right).
\end{eqnarray*}
Using (\ref{odinh22}), we infer from the last inequality that
\begin{eqnarray}\label{2-12-15ct1}
C\left\|\mathcal{U}_{\Gamma_0} -
\mathcal{U}_{\Gamma^\dag}\right\|^2_{H^1(\Omega)}
&\le& \limn \liminf_{m\to\infty} {\mathcal{J}^{h_n}_{\delta_n}}' \left(\Gamma_{n_m}\right) \left( \Gamma_{n_m} - \Gamma_n\right).
\end{eqnarray}
In view of (\ref{23-10-15ct1}) we get that
\begin{eqnarray}\label{2-12-15ct2}
{\mathcal{J}^{h_n}_{\delta_n}}' \left(\Gamma_{n_m}\right) \left( \Gamma_{n_m} - \Gamma_n\right)
&=& \int_\Omega \left( Q_{n_m} - Q_n\right) \nabla \bar{\Pi}^{h_n} z_{\delta_n} \cdot \nabla \bar{\Pi}^{h_n} z_{\delta_n}  \nonumber\\
&~\quad& - 2\left( f_{n_m} - f_n, \bar{\Pi}^{h_n} z_{\delta_n} \right) - 2 \left\langle g_{n_m} - g_n, \gamma\bar{\Pi}^{h_n} z_{\delta_n}\right\rangle \nonumber\\ 
&~\quad& -\int_\Omega \left( Q_{n_m} - Q_n\right) \nabla \mathcal{U}^{h_n}_{\Gamma_{n_m}} \cdot \nabla \mathcal{U}^{h_n}_{\Gamma_{n_m}} + 2\left( f_{n_m} - f_n,  \mathcal{U}^{h_n}_{\Gamma_{n_m}}\right) + 2\left\langle g_{n_m} - g_n,  \gamma \mathcal{U}^{h_n}_{\Gamma_{n_m}}\right\rangle \nonumber\\
&:=&  A_1 -2A_2-2A_3-A_4+2A_5+2A_6. 
\end{eqnarray}
Since $Q_{n_m}\rightharpoonup Q_0$ weakly$^*$ in $\mathbf{L}^\infty_{ \mbox{\tiny sym}}(\Omega)$ as $m\to\infty$, we have for the first term that
\begin{eqnarray*}
\limn \lim_{m\to\infty} A_1
&:= & \limn \left( \lim_{m\to\infty} \int_\Omega \left( Q_{n_m} - Q_n\right) \nabla \bar{\Pi}^{h_n} z_{\delta_n} \cdot \nabla \bar{\Pi}^{h_n} z_{\delta_n}\right) \\
&=& \limn \int_\Omega \left( Q_0 - Q_n\right) \nabla \bar{\Pi}^{h_n} z_{\delta_n} \cdot \nabla \bar{\Pi}^{h_n} z_{\delta_n}\\
&=& \limn \int_\Omega \left( Q_0 - Q_n\right) \nabla \mathcal{U}_{\Gamma^\dag} \cdot \nabla \mathcal{U}_{\Gamma^\dag}\\
&~\quad& + \limn \int_\Omega \left( Q_0 - Q_n\right) \nabla \left( \bar{\Pi}^{h_n} z_{\delta_n} - \mathcal{U}_{\Gamma^\dag}\right) \cdot \nabla \left(\bar{\Pi}^{h_n} z_{\delta_n} + \mathcal{U}_{\Gamma^\dag}\right)\\
&=& \limn \int_\Omega \left( Q_0 - Q_n\right) \nabla \left( \bar{\Pi}^{h_n} z_{\delta_n} - \mathcal{U}_{\Gamma^\dag}\right) \cdot \nabla \left(\bar{\Pi}^{h_n} z_{\delta_n} + \mathcal{U}_{\Gamma^\dag}\right),
\end{eqnarray*}
since $\limn \int_\Omega \left( Q_0 - Q_n\right) \nabla \mathcal{U}_{\Gamma^\dag} \cdot \nabla \mathcal{U}_{\Gamma^\dag} =0$, due to $Q_{n}\rightharpoonup Q_0$ weakly$^*$ in $\mathbf{L}^\infty_{ \mbox{\tiny sym}}(\Omega)$. Furthermore, by (\ref{2-12-15ct4}), we get that
\begin{align*}
&\limn \left| \int_\Omega \left( Q_0 - Q_n\right) \nabla \left( \bar{\Pi}^{h_n} z_{\delta_n} - \mathcal{U}_{\Gamma^\dag}\right) \cdot \nabla \left(\bar{\Pi}^{h_n} z_{\delta_n} + \mathcal{U}_{\Gamma^\dag}\right) \right|  \\
&~\quad \le\quad \limn C\left\| \nabla\big(\bar{\Pi}^{h_n} z_{\delta_n} - \mathcal{U}_{\Gamma^\dag}\big)\right\|_{L^2(\Omega)} \quad=\quad \limn C\left\| \nabla\big(\Pi^{h_n} z_{\delta_n} - \mathcal{U}_{\Gamma^\dag}\big)\right\|_{L^2(\Omega)} \\
&~\quad \le\quad  C\limn \left\| \Pi^{h_n} z_{\delta_n} - \mathcal{U}_{\Gamma^\dag}\right\|_{H^1(\Omega)} \quad=\quad 0.
\end{align*}
Therefore,
\begin{eqnarray}\label{2-12-15ct3}
\limn \lim_{m\to\infty} A_1 \= 0.
\end{eqnarray}
On the other hand, we get 
\begin{eqnarray}\label{2-12-15ct5}
\limn \lim_{m\to\infty} A_2 
&:=& \limn \lim_{m\to\infty}\left( f_{n_m} - f_n, \bar{\Pi}^{h_n} z_{\delta_n} \right) \quad=\quad \limn \left( f_0 - f_n, \bar{\Pi}^{h_n} z_{\delta_n} \right) \nonumber\\
&=& \underbrace{\limn \left( f_0 - f_n, \mathcal{U}_{\Gamma^\dag}\right)}_{=0} + \limn \left( f_0 - f_n, \bar{\Pi}^{h_n} z_{\delta_n} - \mathcal{U}_{\Gamma^\dag}\right) \nonumber\\
&\le& C\limn \left\| \bar{\Pi}^{h_n} z_{\delta_n} - \mathcal{U}_{\Gamma^\dag}\right\|_{L^2(\Omega)}\nonumber\\
&\le& C\limn \left\| \nabla\big(\bar{\Pi}^{h_n} z_{\delta_n} - \mathcal{U}_{\Gamma^\dag}\big)\right\|_{L^2(\Omega)} \quad=\quad 0.
\end{eqnarray}
We now have that
\begin{eqnarray*}
\limn \lim_{m\to\infty} A_3 &:=& \limn \lim_{m\to\infty} \left\langle g_{n_m} - g_n, \gamma\bar{\Pi}^{h_n} z_{\delta_n}\right\rangle \quad=\quad \limn \left\langle g_0 - g_n, \gamma\bar{\Pi}^{h_n} z_{\delta_n}\right\rangle \\
&=& \limn \left\langle g_0 - g_n, \gamma\Pi^{h_n} z_{\delta_n}\right\rangle - |\partial\Omega|^{-1}\limn \left\langle g_0 - g_n, \left\langle 1, \gamma\Pi^{h_n} z_{\delta_n}\right\rangle \right\rangle
\end{eqnarray*}
with
\begin{align*}
&\limn \left\langle g_0 - g_n, \gamma\Pi^{h_n} z_{\delta_n}\right\rangle \\
&~\quad =\quad \limn \left\langle g_0 - g_n, \gamma \big(\Pi^{h_n} z_{\delta_n} - \mathcal{U}_{\Gamma^\dag}\big)\right\rangle + \underbrace{\limn \left\langle g_0 - g_n, \gamma\mathcal{U}_{\Gamma^\dag}\right\rangle}_{=0}\\
&~\quad \le\quad C\limn\|g_0-g_n\|_{L^2(\partial\Omega)}\left\|\gamma\right\|_{\mathcal{L}\big(H^1(\Omega),H^{1/2}(\partial\Omega)\big)}\left\| \Pi^{h_n} z_{\delta_n} - \mathcal{U}_{\Gamma^\dag}\right\|_{H^1(\Omega)}\\
&~\quad \le\quad C\limn\left\| \Pi^{h_n} z_{\delta_n} - \mathcal{U}_{\Gamma^\dag}\right\|_{H^1(\Omega)} \quad=\quad 0
\end{align*}
and
\begin{eqnarray*}
\limn \left\langle g_0 - g_n, \left\langle 1, \gamma\Pi^{h_n} z_{\delta_n}\right\rangle \right\rangle 
&\le& \limn \left|\left\langle 1,\gamma\Pi^{h_n} z_{\delta_n}\right\rangle \right| \left|\left\langle g_0 - g_n, 1\right\rangle\right| \nonumber\\
&\le& C\limn \left\| \Pi^{h_n} z_{\delta_n}\right\|_{H^1(\Omega)}\left|\left\langle g_0 - g_n, 1\right\rangle\right|\\
&\le& C\limn\left|\left\langle g_0 - g_n, 1\right\rangle\right| \quad=\quad 0
\end{eqnarray*}
so that
\begin{eqnarray}\label{2-12-15ct6}
\limn \lim_{m\to\infty} A_3 \= 0.
\end{eqnarray}
Next, we rewrite
\begin{eqnarray*}
\limn \lim_{m\to\infty} A_4 
&:=& \limn \lim_{m\to\infty} \int_\Omega \left( Q_{n_m} - Q_n\right) \nabla \mathcal{U}^{h_n}_{\Gamma_{n_m}} \cdot \nabla \mathcal{U}^{h_n}_{\Gamma_{n_m}} \nonumber\\
&=& \limn \lim_{m\to\infty} \int_\Omega \left( Q_{n_m} - Q_n\right) \nabla \mathcal{U}^{h_n}_{\Gamma_0} \cdot \nabla \mathcal{U}^{h_n}_{\Gamma_0} \nonumber\\
&~\quad& + \limn \lim_{m\to\infty} \int_\Omega \left( Q_{n_m} - Q_n\right) \nabla \left( \mathcal{U}^{h_n}_{\Gamma_{n_m}}- \mathcal{U}^{h_n}_{\Gamma_0} \right) \cdot \nabla \left( \mathcal{U}^{h_n}_{\Gamma_{n_m}} + \mathcal{U}^{h_n}_{\Gamma_0} \right).
\end{eqnarray*}
By (\ref{24-10-15ct1}), likewise as (\ref{2-12-15ct3}), we get that
\begin{eqnarray*}
\limn \lim_{m\to\infty} \int_\Omega \left( Q_{n_m} - Q_n\right) \nabla \mathcal{U}^{h_n}_{\Gamma_0} \cdot \nabla \mathcal{U}^{h_n}_{\Gamma_0} \quad=\quad 0.
\end{eqnarray*}
Furthermore, we have
\begin{eqnarray*}
\left| \int_\Omega \left( Q_{n_m} - Q_n\right) \nabla \left( \mathcal{U}^{h_n}_{\Gamma_{n_m}}- \mathcal{U}^{h_n}_{\Gamma_0} \right) \cdot \nabla \left( \mathcal{U}^{h_n}_{\Gamma_{n_m}} + \mathcal{U}^{h_n}_{\Gamma_0} \right)\right| &\le& C \left\| \mathcal{U}^{h_n}_{\Gamma_{n_m}}- \mathcal{U}^{h_n}_{\Gamma_0} \right\|_{H^1(\Omega)}.
\end{eqnarray*}
By Lemma \ref{finite-dim}, for each fixed $n$ we have that the sequence $\left( \mathcal{U}^{h_n}_{\Gamma_{n_m}}\right)_m \subset  \mathcal{V}^{h_n}_{1}$ converges to $\mathcal{U}^{h_n}_{\Gamma_0}$ in the $H^1(\Omega)$-norm as $m$ tends to $\infty$. Then we deduce that
\begin{align*}
&\limn \lim_{m\to\infty} \left| \int_\Omega \left( Q_{n_m} - Q_n\right) \nabla \left( \mathcal{U}^{h_n}_{\Gamma_{n_m}}- \mathcal{U}^{h_n}_{\Gamma_0} \right) \cdot \nabla \left( \mathcal{U}^{h_n}_{\Gamma_{n_m}} + \mathcal{U}^{h_n}_{\Gamma_0} \right)\right| \\
&~\quad \le\quad C\limn \lim_{m\to\infty} \left\| \mathcal{U}^{h_n}_{\Gamma_{n_m}}- \mathcal{U}^{h_n}_{\Gamma_0} \right\|_{H^1(\Omega)} \quad=\quad C\limn \left\| \mathcal{U}^{h_n}_{\Gamma_0}- \mathcal{U}^{h_n}_{\Gamma_0} \right\|_{H^1(\Omega)} \quad=\quad 0.
\end{align*}
Thus, we obtain
\begin{eqnarray}\label{2-12-15ct8}
\limn \lim_{m\to\infty} A_4 \= 0.
\end{eqnarray}
Finally, we also get that
\begin{eqnarray}\label{2-12-15ct9}
\limn \lim_{m\to\infty} A_5 
&:=& \limn \lim_{m\to\infty} \left( f_{n_m} - f_n, \mathcal{U}^{h_n}_{\Gamma_{n_m}}\right) \nonumber\\
&=& \limn \lim_{m\to\infty} \left( f_{n_m} - f_n, \mathcal{U}^{h_n}_{\Gamma_0}\right) + \limn \lim_{m\to\infty} \left( f_{n_m} - f_n, \mathcal{U}^{h_n}_{\Gamma_{n_m}} - \mathcal{U}^{h_n}_{\Gamma_0}\right) \nonumber\\
&\le& \limn \left( f_0 - f_n, \mathcal{U}^{h_n}_{\Gamma_0}\right) + C\limn \lim_{m\to\infty} \left\| \mathcal{U}^{h_n}_{\Gamma_{n_m}}- \mathcal{U}^{h_n}_{\Gamma_0} \right\|_{H^1(\Omega)} \quad=\quad 0
\end{eqnarray}
and
\begin{eqnarray}\label{2-12-15ct10}
\limn \lim_{m\to\infty} A_6 
&:=& \limn \lim_{m\to\infty}  \left\langle g_{n_m} - g_n,  \gamma \mathcal{U}^{h_n}_{\Gamma_{n_m}}\right\rangle \nonumber\\
&=& \limn \lim_{m\to\infty}  \left\langle g_{n_m} - g_n,  \gamma \big( \mathcal{U}^{h_n}_{\Gamma_{n_m}} - \mathcal{U}^{h_n}_{\Gamma_0}\big)\right\rangle \nonumber\\
&\le& C\limn \lim_{m\to\infty} \left\| \gamma \big( \mathcal{U}^{h_n}_{\Gamma_{n_m}} - \mathcal{U}^{h_n}_{\Gamma_0}\big)\right\|_{L^2(\partial\Omega)} \nonumber\\
&\le& C\limn \lim_{m\to\infty} \left\| \mathcal{U}^{h_n}_{\Gamma_{n_m}}- \mathcal{U}^{h_n}_{\Gamma_0} \right\|_{H^1(\Omega)} \quad=\quad 0.
\end{eqnarray}
Therefore, it follows from the equations (\ref{2-12-15ct2})--(\ref{2-12-15ct10}) that
$$\limn \lim_{m\to\infty} {\mathcal{J}^{h_n}_{\delta_n}}' \left(\Gamma_{n_m}\right) \left( \Gamma_{n_m} - \Gamma_n\right) \quad=\quad 0.$$
Combining this with (\ref{2-12-15ct1}), we obtain that
$
\mathcal{U}_{\Gamma_0} = \mathcal{U}_{\Gamma^\dag}.
$
Then, by the definition of $\Gamma^\dag$, the weakly l.s.c. property of $\mathcal{R}$ and (\ref{odinh22*}), we get
\begin{eqnarray*}
\mathcal{R}(\Gamma^\dag) \quad\le\quad \mathcal{R}(\Gamma_0) \quad\le\quad \liminf_n \mathcal{R}\left(\Gamma_n\right) \quad\le\quad \limsup_n \mathcal{R}\left(\Gamma_n\right) \quad\le\quad \mathcal{R}(\Gamma^\dag).
\end{eqnarray*}
Thus,
$
\mathcal{R}(\Gamma^\dag) = \mathcal{R}(\Gamma_0) = \limn \mathcal{R}\left(\Gamma_n\right).
$
By the uniqueness of $\Gamma^\dag$, we have $\Gamma_0 = \Gamma^\dag$. Furthermore, since $\left(\Gamma_n\right)_n$
weakly converges in $\mathbf{L}^2_{ \mbox{\tiny sym}}(\Omega) \times L^2(\Omega)\times L^2(\partial\Omega)$ to $\Gamma_0$, we conclude from the last equation that $\left(\Gamma_n\right)_n$
converges to $\Gamma_0$ in the $\mathbf{L}^2_{ \mbox{\tiny sym}}(\Omega) \times L^2(\Omega)\times L^2(\partial\Omega)$-norm. 

It remains to show that the sequence $\big(\mathcal{U}^{h_n}_{\Gamma_n}\big)_n$ converges to $\Phi^\dag = \mathcal{U}_{\Gamma^\dag}$ in the $H^1(\Omega)$-norm. We first get from (\ref{24-10-15ct1}) that
\begin{eqnarray}\label{30-9-16ct2}
\limn \big \| \mathcal{U}_{\Gamma^\dag} - \mathcal{U}^{h_n}_{\Gamma^\dag} \big\|_{H^1(\Omega)} \= 0.
\end{eqnarray}
Furthermore, in view of (\ref{30-9-16ct3}) we also have that
\begin{eqnarray}\label{30-9-16ct4}
\frac{C_\Omega \underline{q}}{1+ C_\Omega}\big \| \mathcal{U}^{h_n}_{\Gamma_n} - \mathcal{U}^{h_n}_{\Gamma^\dag} \big\|^2_{H^1(\Omega)} 
&\le& \int_\Omega (Q^\dag-Q_n)\nabla \mathcal{U}^{h_n}_{\Gamma^\dag} \cdot \nabla \big( \mathcal{U}^{h_n}_{\Gamma_n} - \mathcal{U}^{h_n}_{\Gamma^\dag} \big) \\ 
&~\quad& +\big( f_n-f^\dag, \mathcal{U}^{h_n}_{\Gamma_n} - \mathcal{U}^{h_n}_{\Gamma^\dag}\big) + \left\langle g_n-g^\dag, \gamma\big( \mathcal{U}^{h_n}_{\Gamma_n} - \mathcal{U}^{h_n}_{\Gamma^\dag} \big)\right\rangle. \nonumber
\end{eqnarray}
Since $f_n \rightarrow f^\dag$ in the $L^2(\Omega)$-norm and $g_n \rightarrow g^\dag$ in the $L^2(\partial\Omega)$-norm together with the uniform boundedness (\ref{18/5:ct1}), it follows that
\begin{eqnarray}\label{30-9-16ct5}
\limn \left( \big( f_n-f^\dag, \mathcal{U}^{h_n}_{\Gamma_n} - \mathcal{U}^{h_n}_{\Gamma^\dag}\big) + \left\langle g_n-g^\dag, \gamma\big( \mathcal{U}^{h_n}_{\Gamma_n} - \mathcal{U}^{h_n}_{\Gamma^\dag} \big)\right\rangle\right) \=0.
\end{eqnarray}
We now rewrite
\begin{align*}
&\int_\Omega (Q^\dag-Q_n) \nabla \mathcal{U}^{h_n}_{\Gamma^\dag} \cdot \nabla \big( \mathcal{U}^{h_n}_{\Gamma_n} - \mathcal{U}^{h_n}_{\Gamma^\dag} \big)\\
&~\quad =\quad \int_\Omega (Q^\dag-Q_n)\nabla \mathcal{U}_{\Gamma^\dag} \cdot \nabla \big( \mathcal{U}^{h_n}_{\Gamma_n} - \mathcal{U}^{h_n}_{\Gamma^\dag} \big) 
+ \int_\Omega (Q^\dag-Q_n)\nabla \big( \mathcal{U}^{h_n}_{\Gamma^\dag} - \mathcal{U}_{\Gamma^\dag}\big) \cdot \nabla \big( \mathcal{U}^{h_n}_{\Gamma_n} - \mathcal{U}^{h_n}_{\Gamma^\dag} \big).
\end{align*}
We will estimate for two terms in the right hand side of the above equation. For simplicity of notation we here set 
$$Q^\dag-Q_n \quad:=\quad (q^n_{ij})_{i,j=\overline{1,d}}, \quad \nabla \mathcal{U}_{\Gamma^\dag} \quad:=\quad (U_1,...,U_d) \quad\mbox{and}\quad \nabla \big(\mathcal{U}^{h_n}_{\Gamma_n} - \mathcal{U}^{h_n}_{\Gamma^\dag}\big) \quad:=\quad (V^n_1,...,V^n_d).$$
Then, we have
\begin{align*}
&\int_\Omega (Q^\dag-Q_n)\nabla \mathcal{U}_{\Gamma^\dag} \cdot \nabla \big( \mathcal{U}^{h_n}_{\Gamma_n} - \mathcal{U}^{h_n}_{\Gamma^\dag} \big) 
\quad=\quad \int_\Omega \Big( \sum_{j=1}^d q^n_{1j}U_j, ..., \sum_{j=1}^d q^n_{dj}U_j\Big) \cdot (V^n_1,...,V^n_d)\\
&~\quad\le\quad \left( \int_\Omega \Big( \sum_{j=1}^d q^n_{1j}U_j\Big)^2 + ...+ \Big( \sum_{j=1}^d q^n_{dj}U_j\Big)^2\right)^{1/2} \left( \int_\Omega \Big( V^n_1\Big)^2 + ...+ \Big( V^n_d\Big)^2\right)^{1/2}\\
&~\quad\le\quad \left( \int_\Omega \Big(\sum_{i,j=1}^d (q^n_{ij})^2\Big)\Big( \sum_{i=1}^d U_i^2\Big)\right)^{1/2} \left( \int_\Omega \Big( V^n_1\Big)^2 + ...+ \Big( V^n_d\Big)^2\right)^{1/2}\\
&~\quad=\quad \left( \int_\Omega \|Q^\dag-Q_n\|^2_{\mathcal{S}_d} |\nabla \mathcal{U}_{\Gamma^\dag}|^2\right)^{1/2} \left( \int_\Omega\Big|\nabla \big(\mathcal{U}^{h_n}_{\Gamma_n} - \mathcal{U}^{h_n}_{\Gamma^\dag}\big)\Big|^2\right)^{1/2}\\
&~\quad\le\quad \sqrt{2} \left( \int_\Omega \|Q^\dag-Q_n\|^2_{\mathcal{S}_d} |\nabla \mathcal{U}_{\Gamma^\dag}|^2\right)^{1/2} \left( \int_\Omega\big|\nabla \mathcal{U}^{h_n}_{\Gamma_n} \big|^2 + \int_\Omega\big|\nabla \mathcal{U}^{h_n}_{\Gamma^\dag} \big|^2\right)^{1/2}\\
&~\quad\le\quad C\left( \mathcal{R}^2\left(\Gamma_n\right)
+  \mathcal{R}^2\left(\Gamma^\dag\right)\right)^{1/2}\left( \int_\Omega \|Q^\dag-Q_n\|^2_{\mathcal{S}_d} |\nabla \mathcal{U}_{\Gamma^\dag}|^2\right)^{1/2}, \mbox{~by~} (\ref{18/5:ct1})\\
&~\quad\le\quad C\left( \int_\Omega \|Q^\dag-Q_n\|^2_{\mathcal{S}_d} |\nabla \mathcal{U}_{\Gamma^\dag}|^2\right)^{1/2}, \mbox{~by~} (\ref{odinh22*}).
\end{align*}
Similarly, we get
\begin{eqnarray*}
\int_\Omega (Q^\dag-Q_n)\nabla \big( \mathcal{U}^{h_n}_{\Gamma^\dag} - \mathcal{U}_{\Gamma^\dag}\big) \cdot \nabla \big( \mathcal{U}^{h_n}_{\Gamma_n} - \mathcal{U}^{h_n}_{\Gamma^\dag} \big) &\le& C\big \| \mathcal{U}_{\Gamma^\dag} - \mathcal{U}^{h_n}_{\Gamma^\dag} \big\|_{H^1(\Omega)},
\end{eqnarray*}
and arrive at
\begin{align*}
\int_\Omega (Q^\dag-Q_n) \nabla \mathcal{U}^{h_n}_{\Gamma^\dag} \cdot \nabla \big( \mathcal{U}^{h_n}_{\Gamma_n} - \mathcal{U}^{h_n}_{\Gamma^\dag} \big) \quad\le\quad C \left( \int_\Omega \|Q^\dag-Q_n\|^2_{\mathcal{S}_d} |\nabla \mathcal{U}_{\Gamma^\dag}|^2\right)^{1/2} + C\big \| \mathcal{U}_{\Gamma^\dag} - \mathcal{U}^{h_n}_{\Gamma^\dag} \big\|_{H^1(\Omega)}.
\end{align*}
Since $Q_n\rightarrow Q^\dag$ in the $\mathbf{L}^2_{ \mbox{\tiny sym}}(\Omega)$-norm, up to a subsequence
we assume that $(Q_n)_n$  converges to $Q^\dag$ a.e. in $\Omega$. Then, by the Lebesgue dominated convergence theorem, we deduce that
$$\limn \int_\Omega \|Q^\dag-Q_n\|^2_{\mathcal{S}_d} |\nabla \mathcal{U}_{\Gamma^\dag}|^2 \quad=\quad 0.$$
Thus, together with (\ref{30-9-16ct2}), we have
\begin{eqnarray}\label{30-9-16ct6}
\limn \int_\Omega (Q^\dag-Q_n)\nabla \mathcal{U}^{h_n}_{\Gamma^\dag} \cdot \nabla \big( \mathcal{U}^{h_n}_{\Gamma_n} - \mathcal{U}^{h_n}_{\Gamma^\dag} \big) \= 0.
\end{eqnarray}
It follows from (\ref{30-9-16ct4})--(\ref{30-9-16ct6}) that
$\limn \big \| \mathcal{U}^{h_n}_{\Gamma_n} - \mathcal{U}^{h_n}_{\Gamma^\dag} \big\|_{H^1(\Omega)} =0$. By serving of (\ref{30-9-16ct2}) again, we then conclude that $\limn \big \| \mathcal{U}^{h_n}_{\Gamma_n} - \mathcal{U}_{\Gamma^\dag} \big\|_{H^1(\Omega)} =0$, which finishes the proof.
\end{proof}

\section{Error bounds}\label{Convergence rates}

In this section we investigate error bounds of discrete regularized solutions to  the identification problem. For any $\Gamma := (Q,f,g) \in \mathcal{H}_{ad}$ the mapping
$$
\mathcal{U}'_{\Gamma}: \quad\mathbf{L}^\infty_{ \mbox{\tiny sym}}(\Omega) \times
L^{2}(\Omega)\times L^2(\partial\Omega) \quad\rightarrow\quad H^1_\diamond(\Omega)
$$
is linear, continuous with the dual
$$
{\mathcal{U}'_{\Gamma}}^*: \quad {H^1_\diamond(\Omega)}^* \quad\rightarrow\quad
{\mathbf{L}^\infty_{ \mbox{\tiny sym}}(\Omega)}^*\times L^2(\Omega)\times
L^2(\partial\Omega).
$$
\begin{theorem}\label{conv. rate.}
Assume that a function $w^* \in {H^1_\diamond (\Omega)}^*$ exists such that
\begin{eqnarray}\label{moi14***}
{\mathcal{U}'_{\Gamma^\dag}}^*w^* \= \Gamma^\dag.
\end{eqnarray}
Then
\begin{eqnarray} \label{29-6-15ct3}
\left\| \mathcal{U}^h_{\Gamma^h} - \mathcal{U}_{\Gamma^\dag} \right\|^2_{H^1(\Omega)} &+& \rho \left\|\Gamma^h - \Gamma^\dag \right\|^2_{\mathbf{L}^2_{ \mbox{\tiny sym}}(\Omega) \times L^2(\Omega) \times L^2(\partial\Omega)} \nonumber\\
&=& \mathcal{O} \left( h^{-2}\delta^2+ \big(\chi^h_{\Phi^\dag}\big)^2 + \big(\beta^h_{\mathcal{U}_{\Gamma^\dag}}\big)^2 + \left( \chi^h_w\right)^2 +\rho^2\right),
\end{eqnarray}
where $\Gamma^h :=\big(Q^h,f^h,g^h\big)$ is the unique solution to $\left(
\mathcal{P}^{\rho,h}_\delta \right)$ and $w \in H^1_\diamond(\Omega)$ is the unique weak solution of the Neumann problem
\begin{eqnarray}\label{25-9-15ct1}
-\nabla \cdot (Q^\dag\nabla w) \= f^\dag + w^* \quad\mbox{in}\quad \Omega \quad\mbox{and}\quad Q^\dag\nabla w\cdot\vec{n} \quad=\quad g^\dag \quad\mbox{on}\quad \partial\Omega.
\end{eqnarray}
\end{theorem}

\begin{remark}
Due to Remark \ref{15-9-16ct1}, in case $\mathcal{U}_{\Gamma^\dag},~w \in H^2(\Omega)$ we have
$
0\le\chi^h_{\Phi^\dag},~ \beta^h_{\mathcal{U}_{\Gamma^\dag}},~ \chi^h_w \le Ch.
$
Therefore, with $\delta\sim h^2$ and $\rho\sim h$ we obtain the following error bounds
\begin{eqnarray}
&& \left\| \mathcal{U}^h_{\Gamma^h} - \mathcal{U}_{\Gamma^\dag} \right\|_{H^1(\Omega)} \quad=\quad \mathcal{O}(h) \quad\mbox{and}\label{revision0} \\
&&\left\|\Gamma^h - \Gamma^\dag \right\|_{\mathbf{L}^2_{ \mbox{\tiny sym}}(\Omega) \times L^2(\Omega)\times L^2(\partial\Omega)} \quad=\quad \mathcal{O} \big( h^{1/2}\big).\label{revision1}
\end{eqnarray}
\end{remark}

\begin{remark}\label{Goe-1}
Let $\bar{\Gamma} :=(\bar{Q}, \bar{f},\bar{g})\in\mathcal{I}(\Phi^\dag)$ be such that the equation \eqref{moi14***} satisfying with $\bar{\Gamma}$ for some  $w^* \in {H^1_\diamond (\Omega)}^*$, i.e.
$
{\mathcal{U}'_{\bar{\Gamma}}}^*w^* = \bar{\Gamma}.
$ Then $\bar{\Gamma}$ is the unique minimum norm solution of the identification, i.e. $\bar{\Gamma} = \Gamma^\dag$.
\end{remark}
Indeed, due to \eqref{inf1} we have for all $\Gamma := (Q,f,g)\in \mathcal{I}(\Phi^\dag)$ that
\begin{eqnarray*}
\mho \quad:=\quad \left( \bar{\Gamma} , \Gamma -\bar{\Gamma} \right)_{\mathbf{L}^2_{ \mbox{\tiny sym}}(\Omega) \times L^2(\Omega) \times L^2(\partial\Omega)} 
&=& \left\langle \bar{\Gamma} , \Gamma -\bar{\Gamma} \right\rangle_{\left({\mathbf{L}^\infty_{ \mbox{\tiny sym}}(\Omega)}^*\times L^2(\Omega) \times L^2(\partial\Omega), \mathbf{L}^\infty_{ \mbox{\tiny sym}}(\Omega) \times
L^2(\Omega)\times L^2(\partial\Omega)\right)} \nonumber\\
&=& \left\langle {\mathcal{U}'_{\bar{\Gamma}}}^*w^* , \Gamma- \bar{\Gamma} \right\rangle_{\left({\mathbf{L}^\infty_{ \mbox{\tiny sym}}(\Omega)}^*\times L^2(\Omega) \times L^2(\partial\Omega), \mathbf{L}^\infty_{ \mbox{\tiny sym}}(\Omega) \times
L^2(\Omega)\times L^2(\partial\Omega)\right)} \nonumber\\
&=& \left\langle w^*, \mathcal{U}'_{\bar{\Gamma}} \left( \Gamma - \bar{\Gamma}\right)  \right\rangle_{\left({H^1(\Omega)}^*, H^1(\Omega)\right)}\\
&=& \int_\Omega \bar{Q}\nabla \mathcal{U}'_{\bar{\Gamma}} \left( \Gamma - \bar{\Gamma}\right)\cdot \nabla W
\end{eqnarray*}
for some $W\in H^1_\diamond(\Omega)$, since the expression
$ [u,v] := \int_\Omega \bar{Q}\nabla u\cdot \nabla v $
generates a scalar inner product on the space $H^1_\diamond(\Omega)$ which is equivalent to the usual one. By \eqref{22-9-16ct1} we then get 
\begin{eqnarray*}
\mho &=& -\int_\Omega (Q-\bar{Q})\nabla\mathcal{U}_{\bar{\Gamma}}\cdot\nabla W +(f-\bar{f}, W) + \langle g-\bar{g}, \gamma W\rangle\\
&=& \int_\Omega \bar{Q}\nabla\mathcal{U}_{\bar{\Gamma}}\cdot\nabla W -(\bar{f}, W) - \langle \bar{g}, \gamma W\rangle - \left( \int_\Omega Q\nabla\mathcal{U}_{\bar{\Gamma}}\cdot\nabla W -(f, W) - \langle g, \gamma W\rangle \right) \quad=\quad 0, 
\end{eqnarray*}
due to \eqref{ct9} and the fact $\mathcal{U}_{\bar{\Gamma}} = \mathcal{U}_{\Gamma} = \Phi^\dag$. Therefore, we deduce that
\begin{eqnarray*}
\frac{1}{2}\left\| \Gamma  \right\|^2_{\mathbf{L}^2_{ \mbox{\tiny sym}}(\Omega) \times L^2(\Omega) \times L^2(\partial\Omega)} &-& \frac{1}{2}\left\| \bar{\Gamma}  \right\|^2_{\mathbf{L}^2_{ \mbox{\tiny sym}}(\Omega) \times L^2(\Omega) \times L^2(\partial\Omega)} \\
&=& \frac{1}{2}\left\| \Gamma -\bar{\Gamma}  \right\|^2_{\mathbf{L}^2_{ \mbox{\tiny sym}}(\Omega) \times L^2(\Omega) \times L^2(\partial\Omega)} + \mho \quad\ge\quad 0,
\end{eqnarray*}
which completed the proof.

\begin{proof}[Proof of Theorem \ref{conv. rate.}]
Due to the optimality of $\Gamma^h$, we get that
\begin{eqnarray*}
\mathcal{J}^h_\delta \left(\Gamma^h\right) + \rho \mathcal{R}\left(\Gamma^h\right)
&\le& \mathcal{J}^h_\delta\left(\Gamma^\dag\right)
+ \rho \mathcal{R}\left(\Gamma^\dag\right)
\end{eqnarray*}
which implies
\begin{eqnarray} \label{25-10-15ct1}
&&\mathcal{J}^h_\delta \left(\Gamma^h\right) + \rho \left\|\Gamma^h - \Gamma^\dag \right\|^2_{\mathbf{L}^2_{ \mbox{\tiny sym}}(\Omega) \times L^2(\Omega) \times L^2(\partial\Omega)} \nonumber\\
&&~\quad  \le\quad \mathcal{J}^h_\delta\left(\Gamma^\dag\right)
+ 2\rho \left\langle \Gamma^\dag , \Gamma^\dag - \Gamma^h \right\rangle_{\mathbf{L}^2_{ \mbox{\tiny sym}}(\Omega) \times L^2(\Omega) \times L^2(\partial\Omega)}\nonumber\\
&&~\quad \le\quad C \left( h^{-2}\delta^2+ \big(\chi^h_{\Phi^\dag}\big)^2 + \big(\beta^h_{\mathcal{U}_{\Gamma^\dag}}\big)^2\right) + 2\rho \left( \Gamma^\dag , \Gamma^\dag - \Gamma^h \right)_{\mathbf{L}^2_{ \mbox{\tiny sym}}(\Omega) \times L^2(\Omega) \times L^2(\partial\Omega)},
\end{eqnarray}
by Lemma \ref{8-1-16ct1}. Now, by (\ref{inf1}) and (\ref{moi14***}), we infer that
\begin{eqnarray}\label{25-10-15ct1*}
I &:=& \left( \Gamma^\dag , \Gamma^\dag - \Gamma^h \right)_{\mathbf{L}^2_{ \mbox{\tiny sym}}(\Omega) \times L^2(\Omega) \times L^2(\partial\Omega)} \quad=\quad \left\langle w^*, \mathcal{U}'_{\Gamma^\dag} \left( \Gamma^\dag - \Gamma^h\right)  \right\rangle_{\left({H^1(\Omega)}^*, H^1(\Omega)\right)}.
\end{eqnarray}
Thus, by the definition of the weak solution to (\ref{25-9-15ct1}) and (\ref{22-9-16ct1}), we obtain
\begin{eqnarray*}
I &=& \int_\Omega Q^\dag \nabla \mathcal{U}'_{\Gamma^\dag} \left( \Gamma^\dag - \Gamma^h\right) \cdot \nabla w \underbrace{-\big(f^\dag, \mathcal{U}'_{\Gamma^\dag} \left( \Gamma^\dag - \Gamma^h\right)\big) - \big\langle g^\dag, \gamma \mathcal{U}'_{\Gamma^\dag} \left( \Gamma^\dag - \Gamma^h\right)\big\rangle}_{-\int_\Omega Q^\dag \nabla \mathcal{U}_{\Gamma^\dag} \cdot \nabla \mathcal{U}'_{\Gamma^\dag} \left( \Gamma^\dag - \Gamma^h\right), \quad\mbox{by}\quad (\ref{ct9}),}\\
&=& \int_\Omega Q^\dag \nabla \mathcal{U}'_{\Gamma^\dag} \left( \Gamma^\dag - \Gamma^h\right) \cdot \nabla \big(w-\mathcal{U}_{\Gamma^\dag}\big)\\
&=& -\int_\Omega \left( Q^\dag - Q^h \right)\nabla \mathcal{U}_{\Gamma^\dag} \cdot \nabla \big(w-\mathcal{U}_{\Gamma^\dag}\big)  + \left( f^\dag - f^h, w-\mathcal{U}_{\Gamma^\dag}\right) + \left\langle g^\dag-g^h,\gamma \big(w-\mathcal{U}_{\Gamma^\dag}\big)\right\rangle \\
&=&\underbrace{- \int_\Omega Q^\dag \nabla \mathcal{U}_{\Gamma^\dag} \cdot \nabla \big(w-\mathcal{U}_{\Gamma^\dag}\big)  + \left( f^\dag,  w-\mathcal{U}_{\Gamma^\dag} \right) + \left\langle g^\dag,\gamma \big(w-\mathcal{U}_{\Gamma^\dag}\big)\right\rangle }_{=\quad 0, \quad\mbox{by}\quad (\ref{ct9}),}  \\
&~\quad& + \int_\Omega Q^h \nabla \mathcal{U}_{\Gamma^\dag} \cdot \nabla \big(w-\mathcal{U}_{\Gamma^\dag}\big) 
\underbrace{- \left( f^h,  w-\mathcal{U}_{\Gamma^\dag}\right) -\left\langle g^h,\gamma \big(w-\mathcal{U}_{\Gamma^\dag}\big)\right\rangle}_{-\int_\Omega  Q^h \nabla \mathcal{U}_{\Gamma^h} \cdot \nabla \big(w-\mathcal{U}_{\Gamma^\dag}\big) } \nonumber\\
&=& \int_\Omega Q^h \nabla \left( \mathcal{U}_{\Gamma^\dag} -\mathcal{U}_{\Gamma^h} \right) \cdot \nabla\big(w-\mathcal{U}_{\Gamma^\dag}\big) 
\end{eqnarray*}
which yields
\begin{eqnarray}\label{25-10-15ct3}
I
&=& \int_\Omega Q^h \nabla \left(\mathcal{U}_{\Gamma^\dag} - \Pi^h z_\delta \right) \cdot \nabla \big(w-\mathcal{U}_{\Gamma^\dag}\big)    + \int_\Omega Q^h \nabla \left( \mathcal{U}^h_{\Gamma^h } -
\mathcal{U}_{\Gamma^h }\right) \cdot \nabla \big(w-\mathcal{U}_{\Gamma^\dag}\big)     \nonumber\\
&~\quad& + \int_\Omega Q^h \nabla \left( \Pi^h z_\delta - \mathcal{U}^h_{\Gamma^h} \right) \cdot \nabla \big(w-\mathcal{U}_{\Gamma^\dag}\big) \quad :=\quad I_1 + I_2 + I_3.
\end{eqnarray}
For $I_1$ we have from (\ref{21-9-16ct3}) that
\begin{eqnarray}\label{25-10-15ct4}
I_1 &:=&  \int_\Omega Q^h \nabla \left(\mathcal{U}_{\Gamma^\dag} - \Pi^h z_\delta \right) \cdot \nabla \big(w-\mathcal{U}_{\Gamma^\dag}\big)  \quad\le\quad C \left\| \mathcal{U}_{\Gamma^\dag} - \Pi^h z_\delta \right\|_{H^1(\Omega)} \quad\le\quad Ch^{-1}\delta+\chi^h_{\Phi^\dag}.
\end{eqnarray}
Due to (\ref{ct9}) and (\ref{10/4:ct1}), we get
$
\int_\Omega Q^h \nabla \left( \mathcal{U}^h_{\Gamma^h } -
\mathcal{U}_{\Gamma^h }\right) \cdot \nabla \Pi^h \big(w-\mathcal{U}_{\Gamma^\dag}\big) 
=0
$ 
and then infer that
\begin{eqnarray}\label{25-10-15ct6}
I_2 &:=&\int_\Omega Q^h \nabla \left( \mathcal{U}^h_{\Gamma^h } -
\mathcal{U}_{\Gamma^h }\right) \cdot \nabla \big(w-\mathcal{U}_{\Gamma^\dag}\big)   \nonumber\\
&=& \int_\Omega Q^h \nabla \left( \mathcal{U}^h_{\Gamma^h } - \mathcal{U}_{\Gamma^h }\right) \cdot \nabla \left( w-\mathcal{U}_{\Gamma^\dag}  - \Pi^h \big(w-\mathcal{U}_{\Gamma^\dag}\big)  \right)  \nonumber\\
&\le& C \left(  \left\| w - \Pi^h w \right\|_{H^1(\Omega)} + \left\| \mathcal{U}_{\Gamma^\dag} - \Pi^h \mathcal{U}_{\Gamma^\dag} \right\|_{H^1(\Omega)} \right) \quad\le\quad C \left( \chi^h_w + \chi^h_{\Phi^\dag}\right) .
\end{eqnarray}
Finally, we have that
\begin{eqnarray}\label{25-10-15ct7}
I_3 &:=& \int_\Omega Q^h \nabla \left( \Pi^h z_\delta - \mathcal{U}^h_{\Gamma^h} \right) \cdot \nabla \big(w-\mathcal{U}_{\Gamma^\dag}\big)  \nonumber\\
&\le& \left( \int_\Omega Q^h \nabla \left( \mathcal{U}^h_{\Gamma^h} - \Pi^h z_\delta \right) \cdot \nabla \left( \mathcal{U}^h_{\Gamma^h} - \Pi^h z_\delta \right) \right)^{1/2} \cdot  \left( \int_\Omega Q^h \nabla \big(w-\mathcal{U}_{\Gamma^\dag}\big) \cdot \nabla \big(w-\mathcal{U}_{\Gamma^\dag}\big) \right)^{1/2} \nonumber\\
&\le& \underbrace{\frac{1}{4\rho} \int_\Omega Q^h \nabla \left( \mathcal{U}^h_{\Gamma^h} - \Pi^h z_\delta \right) \cdot \nabla \left( \mathcal{U}^h_{\Gamma^h} - \Pi^h z_\delta \right) }_{\mathcal{J}^h_\delta\left(\Gamma^h\right)}
+ \rho \int_\Omega Q^h \nabla \big(w-\mathcal{U}_{\Gamma^\dag}\big) \cdot \nabla \big(w-\mathcal{U}_{\Gamma^\dag}\big) \nonumber\\
&\le& \frac{1}{4\rho}\mathcal{J}^h_\delta \left(\Gamma^h\right) + C\rho.
\end{eqnarray}
It follows from (\ref{25-10-15ct3})--(\ref{25-10-15ct7}) that
\begin{eqnarray*}
I &\le& C\left( h^{-1}\delta + \chi^h_{\Phi^\dag} + \chi^h_w + \rho\right) + \frac{1}{4\rho}\mathcal{J}^h_\delta \left(\Gamma^h\right).
\end{eqnarray*}
Thus, together with (\ref{25-10-15ct1})--(\ref{25-10-15ct1*}), we get
\begin{eqnarray*}
\frac{1}{2}\mathcal{J}^h_\delta \left(\Gamma^h\right) 
+ \rho \left\|\Gamma^h - \Gamma^\dag \right\|^2_{\mathbf{L}^2_{ \mbox{\tiny sym}}(\Omega) \times L^2(\Omega) \times L^2(\partial\Omega)} &\le& C \left( h^{-2}\delta^2+ \big(\chi^h_{\Phi^\dag}\big)^2 + \big(\beta^h_{\mathcal{U}_{\Gamma^\dag}}\big)^2 + \left( \chi^h_w\right)^2 +\rho^2\right),
\end{eqnarray*}
which finishes the proof.
\end{proof}

\section{Gradient projection algorithm with Armijo steplength rule}\label{iterative}

In this section we present the
gradient projection algorithm with Armijo steplength rule (cf.\ \cite{kelley,rusz}) for numerical solution of the minimization  problem $\big(\mathcal{P}^{\rho,h}_\delta \big)$.

We first note that for each $\Gamma := (Q,f,g)\in \mathcal{H}_{ad}$, in view of (\ref{23-10-15ct1}), the $\mathcal{L}^2$-gradient of the strictly convex cost function $\Upsilon^{\rho,h}_\delta$ of the problem $\big( \mathcal{P}^{\rho,h}_\delta \big)$ is given by $\nabla \Upsilon^{\rho,h}_\delta (\Gamma):= \big(\Upsilon_Q(\Gamma), \Upsilon_f(\Gamma), \Upsilon_g(\Gamma)\big)$ with
\begin{eqnarray*}
\Upsilon_Q(\Gamma) &=& \nabla \bar{\Pi}^hz_\delta \otimes \nabla \bar{\Pi}^hz_\delta -\nabla \mathcal{U}^h_{\Gamma} \otimes \nabla
\mathcal{U}^h_{\Gamma} +2\rho Q,\\
\Upsilon_f(\Gamma) &=& 2\big(\mathcal{U}^h_{\Gamma} - \bar{\Pi}^hz_\delta + \rho f\big),\\
\Upsilon_g(\Gamma) &=& 2\big(\gamma \big(\mathcal{U}^h_{\Gamma} -\bar{\Pi}^hz_\delta \big) +\rho g\big)
\end{eqnarray*}
and $\bar{\Pi}^h$ generating from $\Pi^h$ according to (\ref{21-9-16ct4}).

The algorithm is then read as: given a step size control $\beta \in (0,1)$, an initial approximation (cf.\ Remark \ref{22-9-16ct2}) $\Gamma_0 := (Q_0,f_0,g_0) \in \mathcal{H}_{ad} \cap \big({\mathcal{V}^h_0}^{d\times d} \times \mathcal{V}^h_1 \times \mathcal{E}^h_1\big)$, number of iteration $N$ and setting $k=0$.
\begin{enumerate}
\item Compute $\mathcal{U}^h_{\Gamma_k} $ from the variational equation
\begin{eqnarray}\label{eqNh}
\int_\Omega Q_k \nabla \mathcal{U}^h_{\Gamma_k} \cdot \nabla \varphi^h \quad=\quad \left( f_k,\varphi^h\right) + \big\langle g_k,\gamma\varphi^h\big\rangle \quad\mbox{for all}\quad \varphi^h\in
\mathcal{V}_{1}^h
\end{eqnarray}
as well as 
\begin{eqnarray}\label{27-3-16ct4}
\Upsilon_{\rho,\delta}^h (\Gamma_k) 
&=& \int_\Omega Q_k\nabla \left(\mathcal{U}^h_{\Gamma_k}-
\bar{\Pi}^h z_\delta\right) \cdot  \nabla \left(\mathcal{U}^h_{\Gamma_k}-
\bar{\Pi}^h z_\delta\right) \notag\\
&~\quad& + \rho\big( \|Q_k\|^2_{\mathbf{L}^2_{ \mbox{\tiny sym}}(\Omega)} + \|f_k\|^2_{L^2(\Omega)} + \|g_k\|^2_{L^2(\partial\Omega)}\big).
\end{eqnarray}
\item Compute the gradient $\nabla \Upsilon^{\rho,h}_\delta (\Gamma_k):= \big(\Upsilon_{Q_k}(\Gamma_k), \Upsilon_{f_k}(\Gamma_k), \Upsilon_{g_k}(\Gamma_k)\big)$ with
\begin{eqnarray*}
\Upsilon_{Q_k}(\Gamma_k) &=& \nabla \bar{\Pi}^hz_\delta \otimes \nabla \bar{\Pi}^hz_\delta -\nabla \mathcal{U}^h_{\Gamma_k} \otimes \nabla
\mathcal{U}^h_{\Gamma_k} +2\rho Q_k,\\
\Upsilon_{f_k}(\Gamma_k) &=& 2\big(\mathcal{U}^h_{\Gamma_k} - \bar{\Pi}^hz_\delta + \rho f_k\big),\\
\Upsilon_{g_k}(\Gamma_k) &=& 2\big(\gamma \big(\mathcal{U}^h_{\Gamma_k} -\bar{\Pi}^hz_\delta \big) +\rho g_k\big).
\end{eqnarray*}
\item Set $\widetilde{\Gamma}_k := \big(\widetilde{Q}_k, \widetilde{f}_k, \widetilde{g}_k\big)$ with $\widetilde{Q}_k (x) := P_{\mathcal{K}} \big( Q_k(x) - \beta\Upsilon_{Q_k}(\Gamma_k)(x)\big)$, $\widetilde{f}_k (x) := f_k(x) - \beta\Upsilon_{f_k}(\Gamma_k)(x)$ and $\widetilde{g}_k (x) := g_k(x) - \beta\Upsilon_{g_k}(\Gamma_k)(x)$.
\begin{enumerate}
\item  Compute
$\mathcal{U}^h_{\widetilde{\Gamma}_k}$ according to (\ref{eqNh}),
$\Upsilon_{\rho,\delta}^h \big(\widetilde{\Gamma}_k\big)$ according to (\ref{27-3-16ct4}), and with $\tau=10^{-4}$
\begin{eqnarray*}
L &:=& \Upsilon_{\rho,\delta}^h \big(\widetilde{\Gamma}_k\big) - \Upsilon_{\rho,\delta}^h (\Gamma_k) + \tau\beta \big( \|\widetilde{Q}_k-Q_k\|^2_{\mathbf{L}^2_{ \mbox{\tiny sym}}(\Omega)} + \|\widetilde{f}_k-f_k\|^2_{L^2(\Omega)} + \|\widetilde{g}_k-g_k\|^2_{L^2(\partial\Omega)}\big).
\end{eqnarray*}
\item If $L\le 0$

\qquad  go to the next step (c) below

else

\qquad  set $\beta := \frac{\beta}{2}$ and then go back (a)
\item Update $\Gamma_k = \widetilde{\Gamma}_k$, set $k=k+1$.
\end{enumerate}
\item  Compute
\begin{eqnarray}\label{7-7-16ct1}
\mbox{Tolerance} &:=& \big\| \nabla \Upsilon^h_{\rho,\delta}(\Gamma_k) \big\|_{\mathbf{L}^2_{ \mbox{\tiny sym}}(\Omega) \times L^2(\Omega)\times L^2(\partial\Omega)}  -\tau_1 -\tau_2\big\| \nabla \Upsilon^h_{\rho,\delta}(\Gamma_0) \big\|_{\mathbf{L}^2_{ \mbox{\tiny sym}}(\Omega) \times L^2(\Omega)\times L^2(\partial\Omega)}\qquad
\end{eqnarray}
with $\tau_1 := 10^{-3}h$ and $\tau_2 := 10^{-2}h$. If $\mbox{Tolerance} \le 0$ or $k>N$, then stop; otherwise go back Step 1.
\end{enumerate}

\section{Numerical implementation}\label{Numerical implement}

For illustrating the theoretical result we consider the Neumann problem
\begin{eqnarray}
-\nabla\cdot \big(Q^\dag\nabla \Phi^\dag\big) &=& f^\dag \mbox{~in~} \Omega, \label{11-1-16ct1}\\
Q^\dag\nabla\Phi^\dag\cdot\vec{n} &=& g^\dag \mbox{~on~} {\partial \Omega} \label{11-1-16ct2} 
\end{eqnarray}
with $\Omega = \{ x = (x_1,x_2) \in {R}^2 ~|~ -1 < x_1, x_2 < 1\}$.

The special constants in the equation (\ref{5/12/12:ct3}) are chosen as $\underline{q} = 0.05$ and $\overline{q}=10$. For discretization
we divide the interval $(-1,1)$ into $\ell$ equal segments, and so the domain $\Omega = (-1,1)^2$ is divided into $2\ell^2$ triangles, where the diameter of each triangle is $h_{\ell} = \frac{\sqrt{8}}{\ell}$. 

We assume that entries of the symmetric diffusion matrix $Q^\dag$ are discontinuous which are defined as
$$q^\dag_{11} \quad:=\quad 2\chi_{\Omega_{11}} + \chi_{\Omega\setminus\Omega_{11}},\quad 
q^\dag_{12} \quad=\quad q^\dag_{21} \quad:=\quad \chi_{\Omega_{12}} \quad\mbox{and}\quad
q^\dag_{22} \quad:=\quad 3\chi_{\Omega_{22}} + 2\chi_{\Omega\setminus\Omega_{22}},$$
where $\chi_D$ is the characteristic functional of the Lebesgue measurable set $D$ and 
\begin{eqnarray*}
\Omega_{11} &:=& \left\{ (x_1, x_2) \in \Omega ~\big|~ |x_1| \le 3/4 \mbox{~and~} |x_2| \le 3/4 \right\},\\
\Omega_{12} &:=& \left\{ (x_1, x_2) \in \Omega ~\big|~ |x_1| + |x_2| \le 3/4 \right\} \quad\mbox{and}\\
\Omega_{22} &:=& \left\{ (x_1, x_2) \in \Omega ~\big|~ x_1^2 + x_2^2 \le  9/16 \right\}.
\end{eqnarray*}
The source functional $f^\dag$ is assumed to be also discontinuous and defined as
$$f^\dag \quad:=\quad \frac{93-2\pi}{48} \chi_{\Omega_{1}} + \frac{45-2\pi}{48}\chi_{\Omega_{2}} -\frac{3+2\pi}{48} \chi_{\Omega \setminus (\Omega_1 \cup \Omega_2)},$$
where
\begin{eqnarray*}
\Omega_1 &:=& \left\{ (x_1, x_2) \in \Omega ~\big|~ 9(x_1 +1/2)^2 + 16(x_2-1/2)^2 \le 1\right\} \quad\mbox{and}\\
\Omega_2 &:=& \left\{ (x_1, x_2) \in \Omega ~\big|~ |x_1-1/2|\le 1/4 \mbox{~and~} |x_2+1/2| \le 1/4\right\}.
\end{eqnarray*}
The Neumann boundary condition $g^\dag$ is chosen with
\begin{eqnarray*}
g^\dag 
&:=&  -2\chi_{[-1,0]\times\{-1\}} + \chi_{(0,1]\times\{-1\}} - \chi_{[-1,0]\times\{1\}} + 2\chi_{(0,1]\times\{1\}}\\
&&~ +3\chi_{\{-1\}\times(-1,0]} - 4\chi_{\{-1\}\times(0,1)} + 4\chi_{\{1\}\times(-1,0]} - 3\chi_{\{1\}\times(0,1)}.
\end{eqnarray*}
The exact state $\Phi^\dag$ is then computed from the finite element equation $KU=F$, where $K$ and $F$ are the stiffness matrix and the load vector associated with the problem (\ref{11-1-16ct1})--(\ref{11-1-16ct2}), respectively. 

We mention that in the above example the sought functions are chosen to be discontinuous. To reconstruct such discontinuous functions one usually employs the total variation regularization which was originally introduced in image denoising by authors of \cite{ROF92}.
This regularization method was proved to be very effective and analyzed by many authors over the last decades for several ill-posed and inverse problems. We also note that the space of all functions with bounded total variation is a {\it non-reflexive} Banach space and the Tikhonov-function of the total variation regularization  is {\it non-differentiable}, which cause some certain difficulties in numerically treating for non-linear, ill-posed inverse problems. In the present work the cost function is convex and differentiable, the convergence history given in Table \ref{b1} and Table \ref{b2} below shows that the algorithm presented in Section \ref{iterative} performs well for the identification problem with the discontinuous coefficents.

We start the computation with the coarsest level $\ell=3$. To this end, for constructing observations with noise of the exact state $\Phi^\dag$ on this coarsest grid we use
$$z_{\delta_\ell} \quad:=\quad \Phi^\dag + \mathcal{N}_{\overline{\delta_\ell}} \quad\mbox{and}\quad \delta_\ell \quad:=\quad \big\|z_{\delta_\ell} - \Phi^\dag\big\|_{L^2(\Omega)},$$
where $\overline{\delta_\ell} = 10\rho_\ell^{1/2} h_\ell^{3/2}$, $\rho_\ell = 10^{-3}h_\ell$ and $\mathcal{N}_{\overline{\delta_\ell}}$ is a $M^{h_\ell}\times 1$-matrix of random numbers in the interval $\big(-\overline{\delta_\ell},\overline{\delta_\ell}\big)$, $M^{h_\ell} =(\ell+1)^2$ is the number of nodes of the triangulation $\mathcal{T}^{h_\ell}$. Therefore, the exact state $\Phi^\dag$ is only measured at 16 nodes of $\mathcal{T}^{h_\ell}$.

We use the algorithm described in \S \ref{iterative} for computing the numerical solution of the problem $\big(\mathcal{P}_{\rho_\ell,\delta_\ell}^{h_\ell} \big)$. The step size control is chosen with $\beta=0.75$.  As the initial approximation we choose
\begin{eqnarray*}
Q_0 &:=& \left[\begin{array}{cc}2&0\\0&2\end{array}\right], \quad f_0 \quad:=\quad \chi_{[-1,0)\times[-1,1]} - \chi_{[0,1]\times[-1,1]} \quad\mbox{and} \\
g_0 &:=& \chi_{[-1,1]\times\{1\}} - \chi_{[-1,1]\times\{-1\}} + \chi_{\{1\}\times(-1,1)} - \chi_{\{-1\}\times(-1,1)}.
\end{eqnarray*}
At each iteration $k$ we compute Tolerance
defined by (\ref{7-7-16ct1}). Then the iteration was stopped if
$\mbox{Tolerance} \le 0$ or the number of iterations reached the maximum iteration count of 800.

After obtaining the numerical solution $\Gamma_\ell=(Q_\ell,f_\ell,g_\ell)$ and the computed numerical state $\mathcal{U}_\ell = \mathcal{U}^{h_{\ell}}_{\Gamma_\ell}$ of the first iteration process with respect to the coarsest level $\ell=3$, we use their interpolations on the next finer mesh $\ell=6$ as an initial approximation  and an observation  of the exact state for the algorithm on this finer mesh, i.e. for the next iteration process with respect to the level $\ell=6$ we employ 
$$(Q_0,f_0,g_0) \quad:=\quad I^{h_6}_1 \Gamma_3 \quad\mbox{and}\quad z_{\delta_6} \quad:=\quad I^{h_6}_1 \mathcal{U}_3 \quad\mbox{with}\quad \delta_6 \quad:=\quad \big\|z_{\delta_6} - \Phi^\dag\big\|_{L^2(\Omega)}$$ 
and $I^{h_\ell}_1$ being the usual node value interpolation operator on $\mathcal{T}^{h_\ell}$, and so on $\ell=12, 24,\ldots$.  We note that the computation process only requires the measurement data of the exact data for the coarsest level $\ell =3$.

The numerical results are summarized in Table \ref{b1} and Table \ref{b2}, where we present the refinement level $\ell$, mesh size $h_\ell$ of the triangulation, regularization parameter $\rho_\ell$, noise $\delta_\ell$ and number of iterates as well as the final $L^2$-error in the coefficients, the final $L^2$ and $H^1$-error in the states, and their experimental order of convergence (EOC), where
$\mbox{EOC}_\Phi := \frac{\ln \Phi(h_1) - \ln \Phi(h_2)}{\ln h_1 - \ln h_2}$
and $\Phi(h)$ is an error function with respect to  $h$. 

\begin{figure}
\caption{Graphs of $\Phi^\dag$, computed numerical state $\mathcal{U}_\ell$ of the algorithm at the 800$^{\mbox{\tiny th}}$ iteration, and the difference to $\Phi^\dag$.}
\label{h1}
\begin{center}
\includegraphics[scale=0.22]{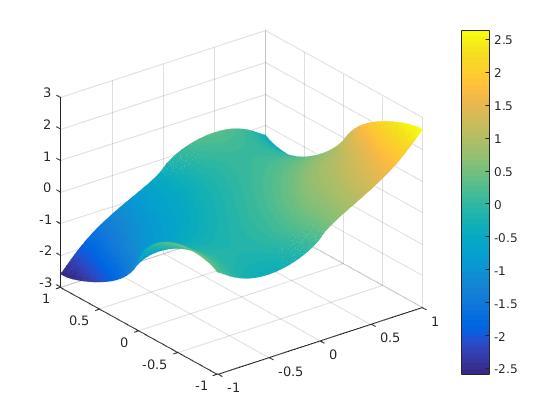}
\includegraphics[scale=0.22]{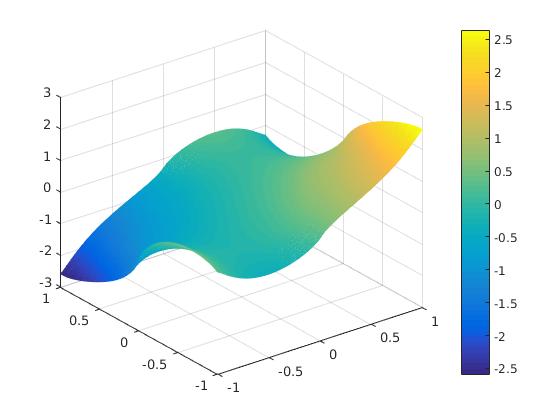}
\includegraphics[scale=0.22]{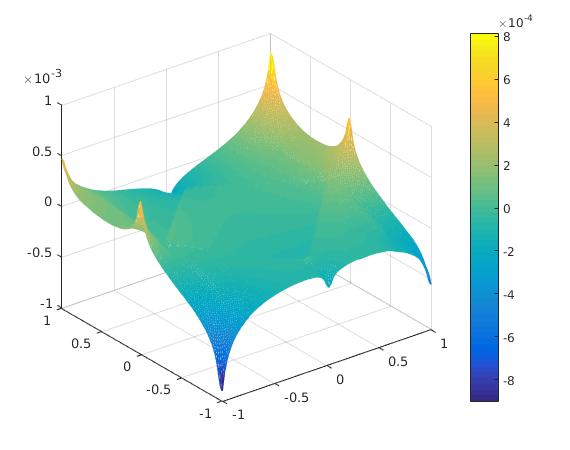}\\
\end{center}
\end{figure}

\begin{figure}
\caption{Graphs of $f_\ell$, $g_\ell$ at the 800$^{\mbox{\tiny th}}$ iteration and the differences $f_\ell - f^\dag$, $g_\ell - g^\dag$.}
\label{h2}
\begin{center}
\includegraphics[scale=0.16]{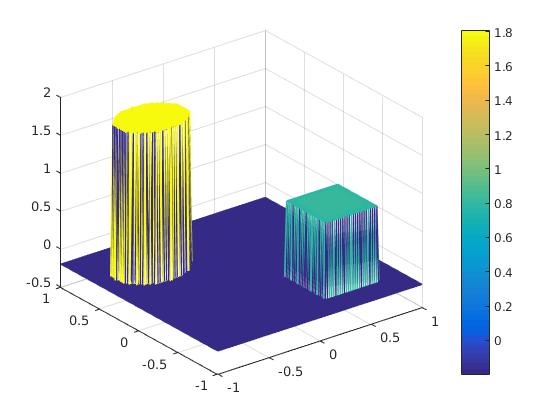}
\includegraphics[scale=0.16]{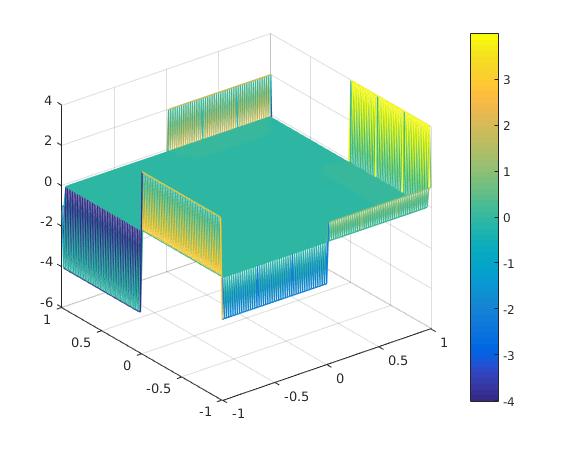}
\includegraphics[scale=0.16]{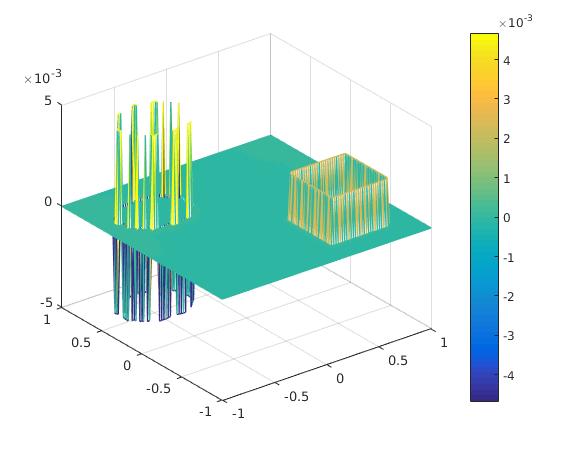}
\includegraphics[scale=0.16]{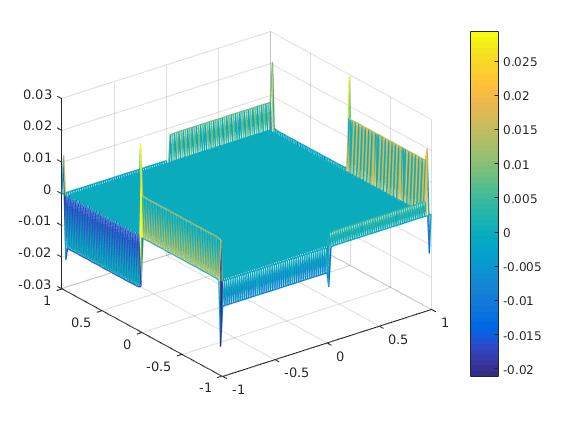}
\end{center}
\end{figure}

\begin{figure}
\caption{Graphs of $q_{\ell,11}$, $q_{\ell,12}$ and $q_{\ell,22}$ at the 800$^{\mbox{\tiny th}}$ iteration.}
\label{h3}
\begin{center}
\includegraphics[scale=0.22]{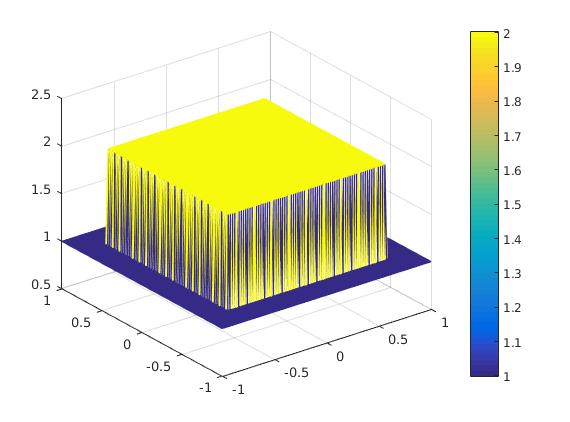}
\includegraphics[scale=0.22]{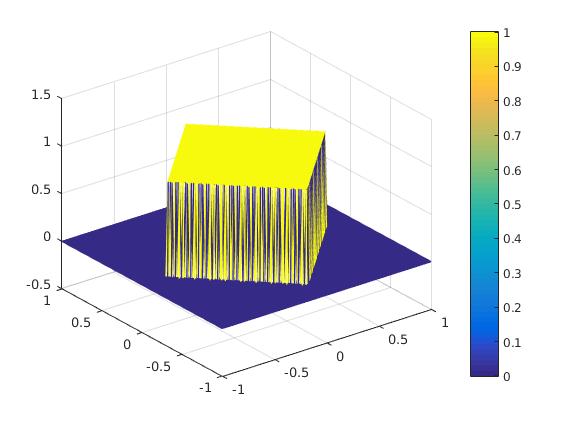}
\includegraphics[scale=0.22]{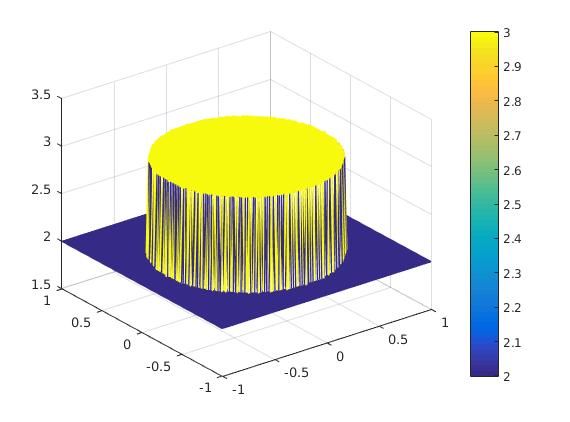}\\
\end{center}
\end{figure}

\begin{figure}
\caption{Differences $q_{\ell,11}-q^\dag_{11}$, $q_{\ell,12}-q^\dag_{12}$ and $q_{\ell,22}-q^\dag_{22}$.}
\label{h4}
\begin{center}
\includegraphics[scale=0.22]{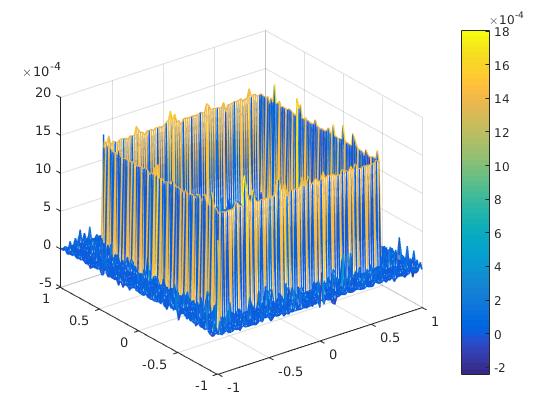}
\includegraphics[scale=0.22]{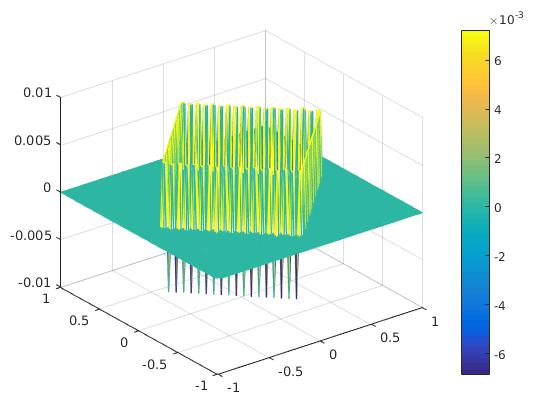}
\includegraphics[scale=0.22]{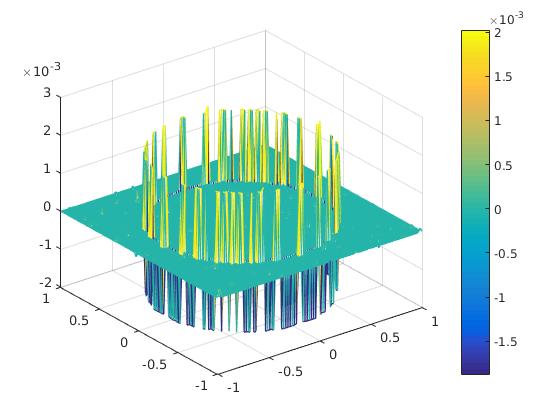}\\
\end{center}
\end{figure}

All figures are here presented corresponding to {\it $\ell = 96$}. Figure \ref{h1} from left to right shows the graphs of $\Phi^\dag$, computed numerical state $\mathcal{U}_\ell$ of the algorithm at the last iteration, and the difference to $\Phi^\dag$. In Figure \ref{h2} we display the computed numerical source term and boundary condition $f_\ell$, $g_\ell$ at the last iteration as well as the differences $f_\ell - f^\dag$, $g_\ell - g^\dag$.  
We write the computed numerical diffusion matrix at the last iteration as
\[Q_\ell \quad:=\quad \left[\begin{array}{cc}q_{\ell,11} & q_{\ell,12}\\
q_{\ell,12} & q_{\ell,22}\end{array}\right].\]
Figure \ref{h3} then shows $q_{\ell,11}$, $q_{\ell,12}$ and $q_{\ell,22}$ while Figure \ref{h4} shows differences $q_{\ell,11}-q^\dag_{11}$, $q_{\ell,12}-q^\dag_{12}$ and $q_{\ell,22}-q^\dag_{22}$. For abbreviation we denote by $\Gamma^\dag :=\big(Q^\dag,f^\dag,g^\dag\big)$ and errors
$$
\Delta \quad:=\quad \left\|\Gamma_\ell - \Gamma^\dag\right\|_{\mathbf{L}^2_{ \mbox{\tiny sym}}(\Omega) \times L^2(\Omega)\times L^2(\partial\Omega)},\quad
\Sigma \quad:=\quad \left\|\mathcal{U}_\ell - \Phi^\dag \right\|_{L^2(\Omega)} \quad\mbox{and}\quad \Lambda \quad:=\quad \left\|\mathcal{U}_\ell - \Phi^\dag \right\|_{H^1(\Omega)}.
$$

\begin{table}
  \caption{Refinement level $\ell$, mesh size $h_\ell$ of the triangulation, regularization parameter $\rho_\ell$, noise $\delta_\ell$ and number of iterates.} \label{b1}
  \begin{center}
    \begin{tabular}{cllll}
    \hline\hline
      $\ell$ & $h_\ell$ & $\rho_\ell$ & 		       $\delta_\ell$ & {\bf Iterate} \\
    \hline
      3   &0.9428    & 9.4281e-4& 0.1755    & 800  \\
\hline
6   &0.4714    & 4.7140e-4& 0.3847    & 800\\
\hline
12  &0.2357    & 2.3570e-4& 0.3334    & 800\\
\hline
24  &0.1179    & 1.1790e-4& 0.1508    & 800\\
\hline
48  &5.8926e-2 & 5.8926e-5& 6.5163e-2 & 800\\
\hline
96  &2.9463e-2 & 2.9463e-5& 2.9896e-2 & 800\\
    \hline\hline
    \end{tabular}
  \end{center}
  \end{table}
  
\begin{table}
  \caption{Errors $\Delta$, $\Sigma$ and $\Lambda$ and Experimental order of convergence between finest and coarsest
level.} \label{b2}
  \begin{center}
    \begin{tabular}{llllll}
    \hline\hline
      {\bf$\Delta$} & {\bf$\Sigma$} & {\bf$\Lambda$} & {\bf EOC$_\Delta$} & {\bf EOC$_\Sigma$} & {\bf EOC$_\Lambda$}\\
    \hline
      0.6349    & 6.2551e-2 & 0.2789    & ---    & ---    & ---   \\
\hline
0.1974    & 3.7602e-2 & 0.1847    & 1.6854 & 0.7342 & 0.5946 \\
\hline
8.3571e-2 & 1.7066e-2 & 0.1382    & 1.2400 & 1.1397 & 0.4184 \\
\hline
3.1600e-2 & 5.4913e-3 & 6.1769e-2 & 1.4031 & 1.6359 & 1.1618 \\
\hline
1.1524e-2 & 9.4491e-4 & 2.0742e-2 & 1.4553 & 2.5389 & 1.5743 \\
\hline
4.1183e-3 & 2.2575e-4 & 8.9372e-3 & 1.4845 & 2.0655 & 1.2147 \\
\hline 
{\bf Mean}&{\bf of}&{ \bf EOC} & 1.4537 & 1.6228 & 0.9928\\
    \hline\hline
    \end{tabular}
  \end{center}
  \end{table}
  
\section*{Acknowledgments} 

The author would like to thank the Referees and the Editor for their valuable comments and suggestions which helped to improve the present paper.


\begin{thebibliography}{9}

\bibitem{Attouch} 
Attouch H., Buttazzo G. and Michaille G., {\em Variational Analysis in Sobolev and BV Space}, Philadelphia: SIAM, 2006.

\bibitem{Baku} 
Baumeister J. and Kunisch K., Identifiability and stability of a two-parameter estimation problem,  {\it Appl.\ Anal.}
\textbf{40}, 263--279, 1991.

\bibitem{banks-kunisch1989} 
Banks H. T. and  Kunisch K.,
{\em Estimation Techniques for Distributed Parameter Systems,
Systems and Control: Foundations and Applications}, Boston: Birkh\"auser, 1989.

\bibitem{Bernardi1} 
Bernardi C., Optimal finite element interpolation on curved domain, {\it SIAM J.\ Numer. Anal.} \textbf{26},  1212--1240, 1989.

\bibitem{Bernardi2} 
Bernardi C. and Girault V., A local regularization operator for triangular and quadrilateral finite elements, {\it SIAM J.\ Numer. Anal.} \textbf{35},  1893--1916, 1998.



\bibitem{Brenner_Scott} 
Brenner S. and Scott R., {\em The Mathematical Theory of Finite Element Methods}, New York: Springer, 2008.

\bibitem{ChanTai2003} 
Chan T.  F and Tai X.  C., Identification of discontinuous coefficients in elliptic problems using total variation regularization, {\it SIAM J.\ Sci.\ Comput.} \textbf{25},  881--904, 2003.

\bibitem{ChanTai2004} 
Chan T.  F and Tai X.  C.,
Level set and total variation regularization for elliptic inverse problems with discontinuous coefficients, {\it J.\ Comput.\ Phys.} \textbf{193}, 40--66, 2004.

\bibitem{chavent2009} 
Chavent G., {\em Nonlinear Least Squares for
Inverse Problems: Theoretical Foundations and Step-by-Step Guide for
Applications}, New York: Springer, 2009.

\bibitem{Chavent_Kunisch2002} 
Chavent G. and Kunisch K., The output least squares identifiability of the diffusion coefficient from an $H^1$-observation in a 2-D elliptic equation, {\it ESAIM Control Optim.\ Calc.\ Var.} \textbf{8},  423--440, 2002.

\bibitem{Chicone} 
Chicone C. and Gerlach J., A note on the identifiability of distributed parameters in elliptic equations, {\it SIAM J.\
Math.\ Anal.} \textbf{18},  1378--1384, 1987.

\bibitem{ciarlet} 
Ciarlet P.  G., {\em Basis Error Estimates for Elliptic Problems}, Handbook of Numerical Analisis, Vol. II,
Ciarlet P. G. and Lions J.-L, eds., Amsterdam:  Elsevier, 1991.

\bibitem{Clement} 
Cl\'ement P., Approximation by finite element functions using local regularization, {\it RAIRO  Anal.\ Num\'er.} \textbf{9}, 77--84, 1975.

\bibitem{Deckelnick} 
Deckelnick K. and Hinze M.,
Convergence and error analysis of a numerical method for the
identification of matrix parameters in elliptic PDEs, {\it Inverse
Problems} \textbf{28}, 15pp, 2012.

\bibitem{Engl_Hanke_Neubauer} 
Engl H.  W., Hanke M. and Neubauer A., {\em Regularization of Inverse Problems:  Mathematics and its Applications}, Dordrecht: Kluwer, 1996.

\bibitem{EnglKuNe} 
Engl H.  W., Kunisch K. and Neubauer A.,
Convergence rates for Tikhonov regularization of nonlinear ill-posed
problems, {\it Inverse Problems} \textbf{5},  523--540, 1989.

\bibitem{Falk} 
Falk R., Error estimates for the numerical identification of a variable coefficient, {\it Math.\ Comput.} \textbf{40}, 537--546, 1983.



\bibitem{hanke}
Hanke M., A regularizing Levenberg-Marquardt scheme, with applications to inverse groundwater filtration problems, {\it Inverse Problems} \textbf{13}, 79--95, 1997.

\bibitem{Haoq}
H{\`a}o D. N. and Quyen T. N. T., Convergence rates for
Tikhonov regularization of coefficient identification problems in
Laplace-type equations, {\it Inverse Problems} \textbf{26}, 23pp, 2010.

\bibitem{hao_quyen3}
H{\`a}o D. N. and Quyen T. N. T., Convergence rates for
Tikhonov regularization of a two-coefficient identification problem
in an elliptic boundary value problem, {\it Numer.\ Math.}
\textbf{120},  45--77, 2012.

\bibitem{hein}
Hein T. and Meyer M., Simultaneous identification of independent parameters in elliptic equations---numerical studies, {\it J.\ Inv.\ Ill Posed Probl.} \textbf{16}, 417--433, 2008.

\bibitem{hinze}
Hinze M., A variational discretization concept in control constrained optimization: the linear-
quadratic case, {\it Comput.\ Optim.\ Appl.} \textbf{30}, 45--61, 2005.

\bibitem{hkq}
Hinze M., Kaltenbacher B. and Quyen T. N. T., Identifying conductivity in electrical impedance tomography with total variation regularization, {\it Numerische Mathematik}, 43pp, 2017 (available https://doi.org/10.1007/s00211-017-0920-8).

\bibitem{Hinze-Tran}
Hinze M. and Quyen T. N. T., Matrix coefficient identification in an elliptic equation with the convex energy functional method. {\it Inverse problems} \textbf{32}, 29pp, 2016.

\bibitem{Hoffmann_Sprekels}
Hoffmann K.  H. and Sprekels J., On the identification of
coefficients of elliptic problems by asymptotic regularization,
{\it Numer.\ Funct.\ Anal.\ Optim.} \textbf{7},  157--177, 1985.

\bibitem{Hsiao}
Hsiao G.  C. and Sprekels J., A stability result for distributed parameter identification in bilinear systems, {\it Math.\ Methods Appl.\ Sci.} \textbf{10},  447--456, 1988.

\bibitem{ito-kunisch}
Ito K.\ and Kunisch K., {\it Lagrange Multiplier Approach to Variational Problems and Applications}, Philadelphia: SIAM, 2008.

\bibitem{bangti2} 
Jin B., Khan T., Maass P.\ and Pidcock M., Function spaces and optimal currents in impedance tomography, {\it J. Inv. Ill-Posed Problems}. \textbf{19}, 25--48, 2011.


\bibitem{isakov}
Isakov V., {\it Inverse Source Problems}, Rhode-Island: American Mathematical Society, 1989.

\bibitem{Kaltenbacher_Schoberl} 
Kaltenbacher B. and Sch\"oberl J.,
A saddle point variational formulation for projection-regularized
parameter identification, {\it Numer.\ Math.} \textbf{91},  675--697, 2002.

\bibitem{kelley}
Kelley C. T., {\it Iterative Methods for Optimization}, Philadelphia.  SIAM, 1999.

\bibitem{keung-zou}
Keung Y. L. and Zou J., An efficient linear solver for nonlinear parameter identification problems, {\it SIAM J.\ Sci.\ Comput} \textbf{22}, 1511--1526, 2000.

\bibitem{knowles}
Knowles I., Uniqueness for an elliptic inverse problem, {\it SIAM J.\ Appl.\ Math.} \textbf{59}, 1356--1370, 1999.

\bibitem{Know3} 
Knowles I. and LaRussa M.  A., Conditional well-posedness for
an elliptic inverse problem, {\it SIAM J.\ Appl.\ Math.} \textbf{71}, 952--971, 2011.

\bibitem{know-wall}
Knowles I. and Wallace R., A variational method for numerical differentiation, {\it Numer.\ Math.} \textbf{70}, 91--110, 1995.

\bibitem{kolo} 
Kohn R.  V. and Lowe B.  D., A variational method for parameter identification, {\it RAIRO Mod\'el.\ Math.\ Anal.\
Num\'er.}  \textbf{22},  119--158, 1988.

\bibitem{Kohn_Vogelius1}
Kohn R.  V. and Vogelius M., Determining conductivity by boundary measurements, {\it Comm.\ Pure Appl.\ Math.} \textbf{37},
289--298, 1984.

\bibitem{Pechstein} 
Pechstein C., {\em Finite and Boundary
Element Tearing and Interconnecting Solvers for Multiscale Problems},
Heidelberg New York Dordrecht London: Springer, 2010.

\bibitem{Rannacher_Vexler} 
Rannacher R. and Vexler B., A priori error estimates for the finite element discretization of elliptic parameter identification problems with pointwise measurements, {\it SIAM J.\ Control Optim.} \textbf{44}, 1844-1863, 2005.

\bibitem{Ric} 
Richter G. R., An inverse problem for the
steady state diffusion equation, {\it SIAM J.\ Appl.\ Math.}
\textbf{41},  210--221, 1981.

\bibitem{ROF92}
Rudin L. I.,  Osher S. J.  and Fatemi E., Nonlinear total variation based noise removal algorithms, {\it Physica} \textbf{D 60}, 259--268, 1992.

\bibitem{rusz}
Ruszczy\'nski A., {\it Nonlinear Optimization}, Princeton: Princeton University Press, 2006.

\bibitem{schuster} 
Schuster T., Kaltenbacher B., Hofmann B. and Kazimierski K. S., {\it Regularization Methods in Banach Spaces}, Berlin: Walter de Gruyter, 2012.

\bibitem{scott_zhang} 
Scott R. and Zhang S. Y., Finite element interpolation of nonsmooth function satisfying boundary conditions, {\it Math.\ Comp.} \textbf{54},  483--493, 1990.

\bibitem{sun}
Sun N.-Z., {\it Inverse Problems in Groundwater Modeling}, Dordrecht: Kluwer, 1994.

\bibitem{Tarantola} 
Tarantola A., {\em Inverse Problem Theory and Methods for Model Parameter Estimation}, Philadelphia: SIAM, 2005.

\bibitem{tartar}
Tartar L., {\it The General Theory of Homogenization}, Berlin: Springer, 2009.

\bibitem{wang_zou} 
Wang L. and Zou J., Error estimates of finite element methods for parameter identification problems in elliptic and parabolic systems, {\it Discrete Contin.\ Dyn.\ Syst.\ Ser.\ B.} \textbf{14},  1641--1670, 2010.

\end{thebibliography}
\end{document}